\documentclass[11pt]{article}
\usepackage{mathrsfs}
\usepackage{amsfonts}
\usepackage{amsmath,amssymb}
\newcommand{\epsln}{\varepsilon}

\newcommand{\diva}{{\rm div}}

\def\bgeq{\begin{equation}}
\def\edeq{\end{equation}}
\def\bgar{\begin{array}}
\def\edar{\end{array}}

\title{Weak Solutions for the Navier-Stokes Equations with
  $B^{-1(\ln)}_{\infty\infty}+B_{X_r}^{-1+r,\frac{2}{1-r}}+L^2$
  Initial Data\thanks{This work is supported by the China national Natural Science
  Foundation under the grant number 11171357.}}
\author{{Shangbin Cui\thanks{E-mail:  cuisb3@yahoo.com.cn}}\\
{\small Department of Mathematics, Sun Yat-Sen University,
Guangzhou, Guandong 510275,}\\ [-0.2cm]
{\small People's Republic of China}\\
}
\date{}

\begin{document}
\maketitle

\begin{abstract}
  In 1934 Leray proved that the Navier-Stokes equations have global weak solutions for
  initial data in $L^2(\mathbb{R}^N)$. In 1990 Calder\'{o}n extended this result to the
  initial value spaces $L^p(\mathbb{R}^N)$ ($2\leq p<\infty$). In the book ``{\em Recent
  developments in the Navier-Stokes problems}'' (2002), Lemari\'{e}-Rieusset extended
  this result of Calder\'{o}n to the space
  $B_{\widetilde{X}_r}^{-1+r,\frac{2}{1-r}}(\mathbb{R}^N)+L^2(\mathbb{R}^N)$ ($0<r<1$),
  where ${X}_r$ is the space of functions whose pointwise products with $H^r$ functions
  belong to $L^2$, $\widetilde{X}_r$ denotes the closure of $C_0^\infty(\mathbb{R}^N)$
  in ${X}_r$, and $B_{\widetilde{X}_r}^{-1+r,\frac{2}{1-r}}(\mathbb{R}^N)$ is the Besov
  space over $\widetilde{X}_r$. In this paper we further extend this result of
  Lemari\'{e}-Rieusset to the larger initial value space $\dot{B}^{-1(\ln)}_{\infty\infty}
  (\mathbb{R}^N)+\dot{B}_{\widetilde{\dot{X}}_r}^{-1+r,\frac{2}{1-r}}(\mathbb{R}^N)+
  L^2(\mathbb{R}^N)$ ($0<r<1$).
\medskip

\medskip

\end{abstract}

\section{Introduction}

\hskip 2em
  In this paper we study existence of solutions for the initial value problem of the
  Navier-Stokes equations:
\begin{eqnarray}
\left\{
\begin{array}{l}
  \partial_t\mathbf{u}-\Delta\mathbf{u}+(\mathbf{u}\cdot\nabla)\mathbf{u}
  +\nabla P=0\quad \mbox{in}\;\,\mathbb{R}^N\times\mathbb{R}_+,\\
  \nabla\cdot\mathbf{u}=0\quad \mbox{in}\;\,\mathbb{R}^N\times\mathbb{R}_+,\\
  \mathbf{u}(x,t)=\mathbf{u}_0(x)\quad \mbox{for}\;\,
  x\in\mathbb{R}^N.
\end{array}
\right.
\end{eqnarray}
%--(1.1)--
  Here $\mathbf{u}=\mathbf{u}(x,t)=(u_1(x,t),u_2(x,t),\cdots,u_N(x,t))$ is an unknown
  $N$-vector function in $(x,t)$ variables, $x\in\mathbb{R}^N$ ($N\geq 2$), $t\geq 0$,
  $P=P(x,t)$ is an unknown scaler function, $\mathbf{u}_0=\mathbf{u}_0(x)$ is a given
  $N$-vector function, $\Delta$ is the Laplacian in the $x$ variables, $\nabla=
  (\partial_{x_1},\partial_{x_2},\cdots,\partial_{x_N})$, and $\mathbb{R}_+=(0,\infty)$
  (later we shall also denote $\overline{\mathbb{R}}_+=[0,\infty)$).

  Existence of solutions of the above problem is a fundamental topic in
  mathematical theory of the Navier-Stokes equations. The first important result on this
  topic was obtained by Leray in 1934 in the reference \cite{LER34}, where he
  introduced the concept of {\em weak solutions} and proved that the problem (1.1) has
  a global weak solution in the class
$$
  C_w([0,\infty),L^2(\mathbb{R}^N))\cap L^\infty([0,\infty),L^2(\mathbb{R}^N))\cap
  L^2([0,\infty),\dot{H}^1(\mathbb{R}^N))
$$
  for any initial data $\mathbf{u}_0\in L^2(\mathbb{R}^N)$, where $C_w([0,\infty),
  L^2(\mathbb{R}^N))$ denotes the set of maps from $[0,\infty)$ to $L^2(\mathbb{R}^N)$
  which are continuous with respect to the weak topology of $L^2(\mathbb{R}^N)$. Here and
  throughout this paper, for simplicity of notations we use the same notation to denote
  a scaler function space and its corresponding $N$-vector counterpart; for instance,
  the notation $L^2(\mathbb{R}^N)$ denotes both the space of scaler $L^2$ functions and
  the space of $N$-vector $L^2$ functions. Since $\mathbf{u}\in L^\infty([0,\infty),
  L^2(\mathbb{R}^N))$, i.e., $\sup_{t\geq 0}\|\mathbf{u}(t)\|_2<\infty$, Leray's weak
  solutions are usually referred to as {\em finite energy} weak solutions in the literature.
  To obtain this result Leray used a smooth approximation approach based on weak
  compactness of bounded sets in separable Banach spaces and dual of Banach spaces. A
  different approach based on Picard's iteration scheme (or Banach's fixed point theorem)
  was introduced by Kato and Fujita in 1964 in \cite{FUJK64}, where they established local
  existence of mild solutions for $\mathbf{u}_0\in H^s(\mathbb{R}^N)$ for $s\geq
  \frac{N}{2}-1$, and global existence of such solutions for small initial dada in
  $H^{\frac{N}{2}-1}(\mathbb{R}^N)$. This approach was later extended to various other
  function spaces, such the Lebesque space $L^p(\mathbb{R}^N)$ for $p\geq N$ by Weissler
  \cite{WEI81}, Kato \cite{KAT84}, Fabes, Johns and Riviere \cite{FABJR72} and Giga
  \cite{GIG86}, the critical and subcritical Sobolev spaces and Besov spaces of either
  positive or negative orders by Kato and Ponce \cite{KATP88},
  Planchon \cite{PLA96}, Terraneo \cite{TER99}, Cannon \cite{CAN95} and et al, the Lorentz
  spaces $L^{p,q}$ by Barraza \cite{BAR96}, the Morrey-Campanato spaces $M^{p,q}$ by
  Giga and Miyakawa \cite{GM89}, the space $BMO^{-1}$ of derivatives of $BMO$ functions
  by Koch and Tataru \cite{KOCT01}, and some general Sobolev and Besov spaces over shift-invariant
  Banach spaces of distributions (see Definition 1.3 below) that can be continuously
  embedded into the Besov space $B_{\infty\infty}^{-1}(\mathbb{R}^N)$ by Lemari\'{e}-Rieusset
  in his expository book \cite{LEM02}. The literatures listed here are far from being
  complete; we refer the reader to see \cite{CAN04} and \cite{LEM02} for expositions and
  references cited therein. The largest initial value space for existence of solutions
  found by using this approach is the $BMO^{-1}$ space introduced by Koch and Tataru in
  \cite{KOCT01}, which coincides with the Triebel-Lizorkin space
  $\dot{F}^{-1}_{\infty 2}(\mathbb{R}^N)$.

  A third approach which combines the above two approaches was introduced by Calder\'{o}n
  in 1990 in \cite{CAL90}. He proved global existence of weak solutions in the class
  $C_w([0,\infty),L^2(\mathbb{R}^N)+L^p(\mathbb{R}^N))$ for initial data $\mathbf{u}_0
  \in L^p(\mathbb{R}^N)$ for any $2\leq p<\infty$. Note that this result of Calder\'{o}n
  fills the gap lying between that of Leray for $p=2$ and those of Weissler et al mentioned
  before for $p\geq N$. In 1998, in an unpublished article
  \cite{LEM98}, in the case $N=3$ Lemari\'{e}-Rieusset extended this result of
  Calder\'{o}n to the initial value space $E_2(\mathbb{R}^3)$, the closure of
  $C_0^\infty(\mathbb{R}^3)$ in the space $L^2_{uloc}(\mathbb{R}^3)$ of uniformly locally
  $L^2$-functions. Later in 2002, in his book \cite{LEM02}, as a concluding result of
  that book, Lemari\'{e}-Rieusset further extended his result to the spaces
  $B_{\widetilde{X}_r}^{s_r,q_r}(\mathbb{R}^N)+L^2(\mathbb{R}^N)$ and
  $B_{\widetilde{X}_r}^{s_r,q_r}(\mathbb{R}^3)+E_2(\mathbb{R}^3)$ ($0<r<1$), where
  $s_r=-1+r$, $q_r=\frac{2}{1-r}$, ${X}_r={X}_r(\mathbb{R}^N)$ is the multiplier space
  of index $r$ (see Definition 1.2 below), $\widetilde{X}_r=\widetilde{X}_r(\mathbb{R}^N)$
  denotes the closure of $C_0^\infty(\mathbb{R}^N)$ in ${X}_r(\mathbb{R}^N)$, and
  $B_{\widetilde{X}_r}^{s_r,q_r}(\mathbb{R}^N)$ is the Besov space over
  $\widetilde{X}_r(\mathbb{R}^N)$.

  On one hand, sine $B_{\widetilde{X}_r}^{s_r,q_r}(\mathbb{R}^N)+L^2(\mathbb{R}^N)
  \not\subseteq bmo^{-1}(\mathbb{R}^N)$, where $bmo^{-1}(\mathbb{R}^N)$ is the
  inhomogeneous version of $BMO^{-1}(\mathbb{R}^N)$, the above-mentioned result of
  Lemari\'{e}-Rieusset gives some new solutions of the Navier-Stokes equations that cannot
  be obtained from either Leray \cite{LER34} or Koch and Tataru \cite{KOCT01}. On the other
  hand, since $r>0$, we see that $s_r=-1+r>-1$, i.e., the regularity index of the space
  $B_{\widetilde{X}_r}^{s_r,q_r}(\mathbb{R}^N)+L^2(\mathbb{R}^N)$ is larger than that of
  $bmo^{-1}(\mathbb{R}^N)={F}^{-1}_{\infty 2}(\mathbb{R}^N)$. Inspired by these observations,
  in this paper we further extend the above result of Lemari\'{e}-Rieusset to the initial
  value space ${B}^{-1(\ln)}_{\infty\infty}(\mathbb{R}^N)+
  {B}_{\widetilde{{X}}_r}^{s_r,q_r}(\mathbb{R}^N)+L^2(\mathbb{R}^N)$, where
  ${B}^{-1(\ln)}_{\infty\infty}(\mathbb{R}^N)$ is a logarithmically modified function space
  to $B^{-1}_{\infty\infty}(\mathbb{R}^N)$ lying in
  between $B^{-1}_{\infty\infty}(\mathbb{R}^N)$ and $B^{s}_{\infty\infty}(\mathbb{R}^N)$
  for any $s>-1$, i.e., $B^{s}_{\infty\infty}(\mathbb{R}^N)\subseteq
  B^{-1(\ln)}_{\infty\infty}(\mathbb{R}^N)\subseteq B^{-1}_{\infty\infty}(\mathbb{R}^N)$
  with continuous embedding (see Definition 1.2 below), and it coincides with the space
  $B^{-1,1}_{\infty\infty}(\mathbb{R}^N)$ introduced by Yoneda in \cite{Yon10}. Note that
  by the recent work of Bourgain and Pavlovi'c \cite{BP08}, the initial value problem for
  the 3-dimensional Navier-Stokes equations is ill-posed in $B^{-1}_{\infty\infty}(\mathbb{R}^3)$
  (see Yoneda \cite{Yon10} for some extensions). In this paper we shall prove that this problem
  is well-posed in ${B}^{-1(\ln)}_{\infty\infty}(\mathbb{R}^N)$; see Theorem 2.1 in Section 2.
  Thus, ${B}^{-1(\ln)}_{\infty\infty}(\mathbb{R}^N)$ is a reasonable substitution of
  $B^{-1}_{\infty\infty}(\mathbb{R}^N)$ in order to get well-posedness of the Navier-Stokes
  initial value problem.

  We point out that since ${B}^{-1(\ln)}_{\infty\infty}(\mathbb{R}^N)\subseteq
  bmo^{-1}(\mathbb{R}^N)$, the result of well-posedness of the Navier-Stokes initial value
  problem in ${B}^{-1(\ln)}_{\infty\infty}(\mathbb{R}^N)$ alone does not provide us with more
  solutions of this problem than those ensured by the result of Koch and Tataru \cite{KOCT01}.
  However, since ${B}^{-1(\ln)}_{\infty\infty}(\mathbb{R}^N)+{B}_{\widetilde{{X}}_r}^{s_r,q_r}
  (\mathbb{R}^N)+L^2(\mathbb{R}^N)\not\subseteq bmo^{-1}(\mathbb{R}^N)$, and this space is
  clearly larger than ${B}_{\widetilde{{X}}_r}^{s_r,q_r}(\mathbb{R}^N)+L^2(\mathbb{R}^N)$,
  our main result essentially enlarges the solution class of the Navier-Stokes initial value
  problem from the known results. The reason that we use the smaller well-posedness space
  ${B}^{-1(\ln)}_{\infty\infty}(\mathbb{R}^N)$ instead of the larger one $bmo^{-1}(\mathbb{R}^N)$
  in our result is because solutions of the Navier-Stokes equations started from
  $bmo^{-1}(\mathbb{R}^N)$ class initial data are not sufficiently regular, so that functions in
  a class of the form $\mathbb{X}+L^2(\mathbb{R}^N)$, where $\mathbb{X}$ represents a general
  function space containing $bmo^{-1}(\mathbb{R}^N)$, cannot be used as an initial value space
  for the Navier-Stokes initial value space. This will be clear from the discussions of Sections
  3 and 4 of the present paper.

  Before presenting the exact statement of our result, let us first make some preparations.
  First we note that since
$$
  (\mathbf{u}\cdot\nabla)\mathbf{u}=\nabla\cdot(\mathbf{u}\otimes\mathbf{u})-
  (\nabla\cdot\mathbf{u})\mathbf{u},
$$
  where $\otimes$ denotes the tenser product between $N$-vectors, the equation in the
  first line of (1.1) can be rewritten as follows:
\begin{eqnarray}
  \partial_t\mathbf{u}-\Delta\mathbf{u}+\nabla\cdot(\mathbf{u}\otimes\mathbf{u})
  +\nabla P=0\quad \mbox{in}\;\,\mathbb{R}^N\times\mathbb{R}_+.
\end{eqnarray}
%--(1.2)--

  {\bf Definition 1.1}\ \ {\em Let $0<T\leq\infty$. Let $\mathbf{u}_0\in
  (\mathcal{D}'(\mathbb{R}^N))^N$ and $\nabla\cdot\mathbf{u}_0=0$. A vector function
  $\mathbf{u}\in(L_{\mbox{\footnotesize loc}}^2(\mathbb{R}^N\times (0,T)))^N\cap C([0,T),
  (\mathcal{D}'(\mathbb{R}^N))^N)$ is said to be a {\em weak solution} of the problem
  $(1.1)$ for $0\leq t<T$, if there exists a distribution $P\in\mathcal{D}'(\mathbb{R}^N
  \times (0,T))$ such that the following two conditions are satisfied:

  $(i)$\ \ $\mathbf{u}$ and $P$ satisfy $(1.2)$ and the equation $\nabla\cdot\mathbf{u}=0$ in
  $\mathbb{R}^N\times (0,T)$ in distribution sense.

  $(ii)$\ $\mathbf{u}(\cdot,0)=\mathbf{u}_0$.}
\medskip

  Let $\mathbb{P}=I+\nabla(-\Delta)^{-1}\nabla$ be the Helmholtz-Weyl projection operator,
  i.e., the $N\times N$ matrix pseudo-differential operator in $\mathbb{R}^N$ with the matrix
  symbol $\Big(\delta_{ij}-\frac{\xi_i\xi_j}{|\xi|^2}\Big)_{i,j=1}^N$, where $\delta_{ij}$
  are the Kronecker symbols. By applying the operator $\mathbb{P}$ to both sides of
  (1.2) we obtain, at least formally, the following equation without the pressure $P$:
\begin{eqnarray*}
  \partial_t\mathbf{u}-\Delta\mathbf{u}+\mathbb{P}\nabla\cdot(\mathbf{u}\otimes\mathbf{u})
  =0\quad \mbox{in}\;\,\mathbb{R}^N\times\mathbb{R}_+.
\end{eqnarray*}
  In order for this equation to make sense, we need to restrict the discussion to {\em
  uniformly locally square integrable weak solutions} (or {\em $L^2_{uloc}$ weak solutions}
  in short), which are weak solutions such that for any $\phi\in C_0^\infty(\mathbb{R}^N
  \times (0,T))$, $\displaystyle\sup_{y\in\mathbb{R}^N}\int_0^T\!\!\int_{\mathbb{R}^N}
  |\phi(x-y,t)\mathbf{u}(x,t)|^2 dxdt<\infty$. For such solutions, the expression
  $\mathbb{P}\nabla\cdot(\mathbf{u}\otimes\mathbf{u})$ is meaningful. Indeed, by choosing
  a function $\varphi\in C_0^\infty(\mathbb{R}^N)$ such that $\varphi(\xi)=1$ for $|\xi|
  \leq 1$ and $\varphi(\xi)=0$ for $|\xi|\geq 2$, we can write
\begin{eqnarray*}
  \mathbb{P}\nabla\cdot(\mathbf{u}\otimes\mathbf{u})
  =\varphi(D)\mathbb{P}\nabla\cdot(\mathbf{u}\otimes\mathbf{u})+
  [I-\varphi(D)]\mathbb{P}\nabla\cdot(\mathbf{u}\otimes\mathbf{u}).
\end{eqnarray*}
  Clearly, $[I-\varphi(D)]\mathbb{P}\nabla$ is a continuous linear map from
  $(\mathcal{S}'(\mathbb{R}^N))^{N\times N}$ to $(\mathcal{S}'(\mathbb{R}^N))^{N}$, so
  that the second term on the right-hand side of the above equality is meaningful. Moreover,
  since each component of $\varphi(D)\mathbb{P}\nabla$ is a convolution operator with a
  $L^1$ kernel, so that it maps $(\mathbf{u}\otimes\mathbf{u})(\cdot,t)\in
  (L^1_{uloc}(\mathbb{R}^N))^{N\times N}$ to $(L^1_{uloc}(\mathbb{R}^N))^{N}$ (for fixed
  $t$), it follows that the first term on the right-hand side also makes sense. Hence,
  the problem (1.1) reduces into the following apparently simpler problem:
\begin{eqnarray}
\left\{
\begin{array}{l}
  \partial_t\mathbf{u}-\Delta\mathbf{u}+\mathbb{P}\nabla\cdot(\mathbf{u}\otimes\mathbf{u})
  =0\quad \mbox{in}\;\,\mathbb{R}^N\times\mathbb{R}_+,\\
  \mathbf{u}(x,t)=\mathbf{u}_0(x)\quad \mbox{for}\;\,
  x\in\mathbb{R}^N.
\end{array}
\right.
\end{eqnarray}
%--(1.3)--
  For equivalence of the problems (1.1) and (1.3) in the category of $L^2_{uloc}$ weak
  solutions, we refer the reader to see Theorem 11.1 of \cite{LEM02}.

  Next we recall the definitions of some function spaces (see \cite{LEM02} for more details):
\begin{itemize}
\item  For $0<r<\frac{N}{2}$, we define the {\em multiplier space of index $r$},
  $X_r(\mathbb{R}^N)$, as the Banach space of locally square-integrable functions on
  $\mathbb{R}^N$ such that pointwise multiplication with these functions maps boundedly
  $H^r(\mathbb{R}^N)$ into $L^2(\mathbb{R}^N)$. The norm of $X_r(\mathbb{R}^N)$ is given
  by
$$
  \|u\|_{X_r}=\sup\{\|uv\|_2:v\in H^r(\mathbb{R}^N),\|v\|_{H^r}\leq 1\}, \quad
  \forall u\in X_r(\mathbb{R}^N).
$$
  $\widetilde{X}_r(\mathbb{R}^N)$ denotes the the closure of $C_0^\infty(\mathbb{R}^N)$ in
  ${X}_r(\mathbb{R}^N)$. We note that $L^\infty(\mathbb{R}^N)\hookrightarrow
  {X}_r(\mathbb{R}^N)\hookrightarrow {X}_s(\mathbb{R}^N)$ for $0<r<s<\frac{N}{2}$.
\item Let $E$ be a shift-invariant Banach space of
  distributions on $\mathbb{R}^N$. For $s\in\mathbb{R}$ and $1\leq q\leq\infty$, we define the
  {\em Besov space over $E$ of index $(s,q)$}, $B^{s,q}_E(\mathbb{R}^N)$, as the Banach space
  of temperate distributions $u$ on $\mathbb{R}^N$ such that $S_0u\in E$, $\Delta_j u\in E$
  for all $j\in\mathbb{N}$, and $\{2^{js}\|\Delta_j u\|_E\}_{j=1}^\infty\in l^q$, where $S_0$
  and $\Delta_j$ are the frequency-localizing operator appearing in the Littlewood-Paley
  decomposition $u=S_0u+\sum_{j=1}^{\infty}\Delta_j u$. The norm of $B^{s,q}_E(\mathbb{R}^N)$
  is given by
$$
  \|u\|_{B^{s,q}_E}=\left\{
\begin{array}{ll}
  \displaystyle\|S_0u\|_E+\Big[\sum_{j=1}^\infty\Big(2^{js}\|\Delta_j u\|_E
  \Big)^q\Big]^{\frac{1}{q}}, \quad &\mbox{if}\;\; 1\leq q<\infty,\\ [0.3cm]
  \displaystyle\|S_0u\|_E+\sup_{j\in\mathbb{N}}\Big(2^{js}\|\Delta_j u\|_E
  \Big), \quad &\mbox{if}\;\; q=\infty,
\end{array}
\right.
   \qquad \forall u\in B^{s,q}_E(\mathbb{R}^N).
$$
  In particular, if $E=L^p(\mathbb{R}^N)$ ($1\leq p\leq\infty$)
  then $B^{s,q}_E(\mathbb{R}^N)=B^s_{pq}(\mathbb{R}^N)$ is the usual Besov space.
\end{itemize}
\medskip

  Let $\gamma\geq 0$ and $\gamma>s$. It is well-known that $\varphi\in
  B^{s,q}_E(\mathbb{R}^N)$ if and only if for any $t>0$ we have $e^{t\Delta}\varphi
  \in E$ and $t^{-\frac{s}{2}+\frac{\gamma}{2}}\|D^{\gamma}e^{t\Delta}\varphi\|_E\in
  L^q((0,t_0),\frac{dt}{t})$, and the norms $\|\varphi\|_{B^{s,q}_E}$ and
$$
  \|\varphi\|_{B^{s,q}_E}':=\|e^{t_0\Delta}\varphi\|_E+
  \Big[\int_0^{t_0}\!\Big(t^{-\frac{s}{2}+\frac{\gamma}{2}}
  \|D^{\gamma}e^{t\Delta}\varphi\|_E\Big)^q\frac{dt}{t}\Big]^{\frac{1}{q}}
$$
  $($for $1\leq q<\infty)$ or $\|\varphi\|_{B^{s,q}_E}':=\|e^{t_0\Delta}\varphi\|_E+
  \sup_{t\in [0,t_0)}t^{-\frac{s}{2}+\frac{\gamma}{2}}\|D^{\gamma}e^{t\Delta}
  \varphi\|_E$ $($for $q=\infty)$ are equivalent (cf. Theorem 5.3 of \cite{LEM02}).
  We thus introduce the following concept:
\medskip

  {\bf Definition 1.2}\ \ {\em Let $E$ be a shift-invariant Banach space of distributions
  on $\mathbb{R}^N$. Let $T>0$. We denote by $B^{-1(\ln),\infty}_{E}$ the function space
  of all distribution $\varphi\in\mathcal{S}'(\mathbb{R}^N)$ such that $e^{t\Delta}\varphi
  \in E$ for all $t\in (0,\infty)$, and
$$
  \|\varphi\|_{B^{-1(\ln),\infty}_{E,T}}:=
  \sup_{t\in (0,T)}\sqrt{t}\Big|\ln\Big(\frac{t}{e^2T}\Big)\Big|
  \|e^{t\Delta}\varphi\|_E<\infty.
$$
  In particular, for $E=L^{\infty}(\mathbb{R}^N)$, the space $B^{-1(\ln),\infty}_{E}$
  is denoted as $B^{-1(\ln)}_{\infty\infty}(\mathbb{R}^N)$. The notations
  $\mathcal{C}_{bu}^{-1(\ln)}(\mathbb{R}^N)$ and $\mathcal{C}_{0}^{-1(\ln)}(\mathbb{R}^N)$
  denote the closures of $C_{bu}(\mathbb{R}^N)$ $($=the space of bounded and uniformly
  continuous functions on $\mathbb{R}^N)$ and  $C_{0}^\infty(\mathbb{R}^N)$ in
  $B^{-1(\ln)}_{\infty\infty}(\mathbb{R}^N)$, respectively. $\quad\Box$}
\medskip

  The definition of $B^{-1(\ln),\infty}_{E}$ does not depend on the specific choice
  of $T$. Indeed, it is easy to check that for any $0<T_1<T_2$ we have
$$
  \|\varphi\|_{B^{-1(\ln),\infty}_{E,T_1}}\leq
  \|\varphi\|_{B^{-1(\ln),\infty}_{E,T_2}}\leq
  \sqrt{\frac{T_2}{T_1}}\|\varphi\|_{B^{-1(\ln),\infty}_{E,T_1}}.
$$
  We note that $B^{-1(\ln),\infty}_{E}$ is a Banach space. Moreover, we have the following
  embedding result:
\medskip

  {\bf Lemma 1.3}\ \ {\em $(i)$\ For any $1<q\leq\infty$, $B^{-1(\ln),\infty}_{E}$
  is continuously embedded into $B^{-1,q}_{E}$.

   $(ii)$\ For any $s>-1$, $B^{s,\infty}_{E}$ is continuously embedded into
   $B^{-1(\ln),\infty}_{E}$.}
\medskip

  {\em Proof}:\ \ For $q=\infty$ the assertion $(i)$ is immediate, because
  $\Big|\ln\Big(\frac{t}{e^2T}\Big)\Big|\geq 2$ for $0<t<T$. For $1<q<\infty$ we
  deduce as follows:
\begin{align*}
  \int_0^{T}\!\Big(t^{\frac{1}{2}}\|e^{t\Delta}\varphi\|_{E}\Big)^q\frac{dt}{t}
  =&\int_0^{T}\!t^{-1}\Big|\ln\Big(\frac{t}{e^2T}\Big)\Big|^{-q}
  \cdot\Big(\sqrt{t}\Big|\ln\Big(\frac{t}{e^2T}\Big)\Big|
  \|e^{t\Delta}\varphi\|_{E}\Big)^qdt
\\
  \leq&\int_0^{T}\!t^{-1}\Big|\ln\Big(\frac{t}{e^2T}\Big)\Big|^{-q}dt
  \cdot\Big(\sup_{0<t<T}\sqrt{t}\Big|\ln\Big(\frac{t}{e^2T}\Big)\Big|
  \|e^{t\Delta}\varphi\|_{E}\Big)^q
\\
  =&\frac{2^{1-q}}{q-1}\Big(\sup_{0<t<T}\sqrt{t}\Big|\ln\Big(\frac{t}{e^2T}\Big)\Big|
  \|e^{t\Delta}\varphi\|_{E}\Big)^q
\end{align*}
  Hence $\|\varphi\|_{B^{-1,q}_{E,T}}\leq(\frac{2^{1-q}}{q-1})^{\frac{1}{q}}
  \|\varphi\|_{B^{-1(\ln),\infty}_{E,T}}$.
  This proves the assertion $(i)$. The assertion $(ii)$ is immediate. $\quad\Box$
\medskip

  {\bf Remark}\ \ It can be easily shown that when $E=L^p(\mathbb{R}^N)$ ($1\leq p\leq\infty$),
  $B^{-1(\ln),\infty}_{E}$ essentially coincides with the space $B^{-1,1}_{p\infty}
  (\mathbb{R}^N)$ introduced by Yoneda in \cite{Yon10}.
\medskip

  The main result of this paper is as follows:
\medskip

  {\bf Theorem 1.4}\ \ {\em Let $0<r<1$, $s_r=-1+r$ and $q_r=\frac{2}{1-r}$. Given $T>0$,
  there exist constants $\epsln_1>0$ and $\epsln_2>0$ such that for any $\mathbf{u}_0=
  \mathbf{v}_0+\mathbf{w}_0+\mathbf{z}_0$ with $\mathbf{v}_0\in{B}^{-1(\ln)}_{\infty\infty}
  (\mathbb{R}^N)$, $\mathbf{w}_0\in{B}_{{X}_r}^{s_r,q_r}(\mathbb{R}^N)$, $\mathbf{z}_0
  \in L^2(\mathbb{R}^N)$, $\nabla\cdot\mathbf{v}_0=\nabla\cdot\mathbf{w}_0=\nabla\cdot
  \mathbf{z}_0=0$, $\|\mathbf{v}_0\|_{{B}^{-1(\ln)}_{\infty\infty,T}}<\epsln_1$ and
  $\|\mathbf{w}_0\|_{{B}_{{X}_r}^{s_r,q_r}}<\epsln_2$, the problem $(1.1)$ has a
  weak solution $\mathbf{u}\in C_w([0,T],{B}^{-1(\ln)}_{\infty\infty}
  (\mathbb{R}^N)+{B}_{{X}_r}^{s_r,q_r}(\mathbb{R}^N)+L^2(\mathbb{R}^N))$. In
  particular, for any $\mathbf{u}_0=\mathbf{v}_0+\mathbf{w}_0+\mathbf{z}_0$ with
  $\mathbf{v}_0\in{\mathcal{C}}^{-1(\ln)}_{0}(\mathbb{R}^N)$, $\mathbf{w}_0\in
  {B}_{\widetilde{{X}}_r}^{s_r,q_r}(\mathbb{R}^N)$, $\mathbf{z}_0\in
  L^2(\mathbb{R}^N)$ and $\nabla\cdot\mathbf{v}_0=\nabla\cdot\mathbf{w}_0=\nabla\cdot
  \mathbf{z}_0=0$,  the problem $(1.1)$ has a weak solution $\mathbf{u}\in
  C_w([0,T],{\mathcal{C}}^{-1(\ln)}_{0}(\mathbb{R}^N)+{B}_{\widetilde{{X}}_r}^{s_r,q_r}
  (\mathbb{R}^N)+L^2(\mathbb{R}^N))$.}
\medskip

  The idea of the proof of this result is as follows. We first solve the
%  which is in the line of Calder\'{o}n \cite{CAL90} and Lemari\'{e}-Rieusset \cite{LEM02},
  Navier-Stokes equations for initial data $\mathbf{v}_0\in{B}^{-1(\ln)}_{\infty\infty}
  (\mathbb{R}^N)$ in a regular space. For this purpose we need to establish well-posedness
  of the Navier-Stokes equations in the space ${B}^{-1(\ln)}_{\infty\infty}(\mathbb{R}^N)$.
  Next we consider the solution of the Navier-Stokes equations of the form $\mathbf{v}
  +\mathbf{w}$ such that $\mathbf{w}$ satisfies the initial condition $\mathbf{w}_0\in
  {B}_{{X}_r}^{s_r,q_r}(\mathbb{R}^N)$. This requires to solve a linear
  perturbation of the Navier-Stokes equations. Generally speaking, such problems are not
  always solvable due to the nonlinearity of the perturbed equations. Fortunately, since
  $\mathbf{v}$ is the solution of the Navier-Stokes equations in a regular space, this
  perturbation problem is solvable and, importantly, the solution also belongs to a
  regular space. Finally we consider the solution of the Navier-Stokes equations of the
  form $\mathbf{v}+\mathbf{w}+\mathbf{z}$ such that $\mathbf{z}$ satisfies the initial
  condition $\mathbf{z}_0\in L^2(\mathbb{R}^N)$. For this purpose we need to solve
  another linear perturbation of the Navier-Stokes equations. Thanks to the fact that
  both $\mathbf{v}$ and $\mathbf{w}$ are in regular spaces, this perturbation problem is
  also solvable. But due to the weak regularity of the initial data $\mathbf{z}_0$,
  this time we can only get a weak solution.

  The organization of the rest part is as follows. In Section 2 we solve the Navier-Stokes
  initial value problem for initial data $\mathbf{v}_0\in{B}^{-1(\ln)}_{\infty\infty}
  (\mathbb{R}^N)$. In Section 3 we solve a linearly perturbed Navier-Stokes initial
  value problem for initial data $\mathbf{w}_0\in{B}_{{X}_r}^{s_r,q_r}
  (\mathbb{R}^N)$. In the last section we solve another linearly perturbed Navier-Stokes
  initial value problem for initial data $\mathbf{z}_0\in L^2(\mathbb{R}^N)$ and give the
  proof of Theorem 1.2.

\section{Well-posedness in $B^{-1(\ln)}_{\infty\infty}(\mathbb{R}^N)$ and
  related spaces}

\hskip 2em

  In this section we prove the following two theorems:
\medskip

  {\bf Theorem 2.1}\ \ {\em The problem $(1.3)$ is semi-globally weakly
  well-posed in $B^{-1(\ln)}_{\infty\infty}(\mathbb{R}^N)$ for small initial data,
  namely, for any $T>0$ there exists corresponding constant $\epsln>0$ such that
  for any $\mathbf{u}_0\in B^{-1(\ln)}_{\infty\infty}(\mathbb{R}^N)$ satisfying the
  conditions $\nabla\cdot\mathbf{u}_0=0$ and $\|\mathbf{u}_0\|_{B^{-1(\ln)}_{\infty\infty}}
  <\epsln$, the problem $(1.1)$ has a unique mild solution in the class
\begin{eqnarray*}
  &\mathbf{u}\in C_w([0,T],B^{-1(\ln)}_{\infty\infty}(\mathbb{R}^N))\cap
  L^\infty_{\rm loc}((0,T),L^\infty(\mathbb{R}^N)), \quad
  \nabla\cdot\mathbf{u}=0,& \\
  &\displaystyle\sup_{t\in (0,T)}\sqrt{t}\Big|\ln\Big(\frac{t}{e^2T}\Big)\Big|
  \|\mathbf{u}(t)\|_{\infty}<\infty,
%  \quad \displaystyle\lim_{t\to 0^+}\sqrt{t}\Big|\ln\Big(\frac{t}{e^2T}\Big)\Big|
%  \|\mathbf{u}(t)\|_{\infty}=0,&
\end{eqnarray*}
  and the solution map $\mathbf{u}_0\mapsto \mathbf{u}$ from $\overline{B}(0,\epsln)
  \subseteq B^{-1(\ln)}_{\infty\infty}(\mathbb{R}^N)$ to the Banach space of the above class
  of functions on $\mathbb{R}^N\times [0,T]$ is Lipschitz continuous. Here $C_{w}([0,T],
  B^{-1(\ln)}_{\infty\infty}(\mathbb{R}^N))$ denotes the Banach space of mappings from
  $[0,T]$ to $B^{-1(\ln)}_{\infty\infty}(\mathbb{R}^N)$ which are continuous with respect to
  the $\ast$-weak topology of $B^{-1(\ln)}_{\infty\infty}(\mathbb{R}^N)$.}
\medskip

  {\bf Theorem 2.2}\ \ {\em $(i)$  The problem $(1.3)$ is locally well-posed in
  $\mathcal{C}^{-1(\ln)}_{bu}(\mathbb{R}^N)$, namely, for any $\mathbf{u}_0\in
  \mathcal{C}^{-1(\ln)}_{bu}(\mathbb{R}^N)$ satisfying the condition $\nabla\cdot\mathbf{u}_0=0$,
  there exists corresponding $T>0$ such that the problem $(1.1)$ has a unique mild solution
  in the class
\begin{eqnarray*}
  &\mathbf{u}\in C([0,T],\mathcal{C}^{-1(\ln)}_{bu}(\mathbb{R}^N))\cap
  L^\infty_{\rm loc}((0,T),L^\infty(\mathbb{R}^N)), \quad
  \nabla\cdot\mathbf{u}=0,& \\
  &\displaystyle\sup_{t\in (0,T)}\sqrt{t}\Big|\ln\Big(\frac{t}{e^2T}\Big)\Big|
  \|\mathbf{u}(t)\|_{\infty}<\infty,
  \quad \displaystyle\lim_{t\to 0^+}\sqrt{t}\Big|\ln\Big(\frac{t}{e^2T}\Big)\Big|
  \|\mathbf{u}(t)\|_{\infty}=0,&
\end{eqnarray*}
  and the solution map $\mathbf{u}_0\mapsto \mathbf{u}$ from $\overline{B}(0,\epsln)
  \subseteq\mathcal{C}^{-1(\ln)}_{bu}(\mathbb{R}^N)$ to the Banach space of the above class of
  functions on $\mathbb{R}^N\times [0,T]$ is Lipschitz continuous.

  $(ii)$  The problem $(1.3)$ is semi-globally well-posed in $\mathcal{C}^{-1(\ln)}_{bu}
  (\mathbb{R}^N)$ for small initial data, namely, for any $T>0$ there exists corresponding
  constant $\epsln>0$ such that for any $\mathbf{u}_0\in\mathcal{C}^{-1(\ln)}_{bu}
  (\mathbb{R}^N)$ satisfying the conditions $\nabla\cdot\mathbf{u}_0=0$ and
  $\|\mathbf{u}_0\|_{B^{-1(\ln)}_{\infty\infty}}<\epsln$, the problem $(1.1)$ has a unique
  mild solution in the above class, and the solution map $\mathbf{u}_0\mapsto\mathbf{u}$
  from $\overline{B}(0,\epsln)\subseteq\mathcal{C}^{-1(\ln)}_{bu}(\mathbb{R}^N)$ to the
  Banach space of the above class of functions on $\mathbb{R}^N\times [0,T]$ is Lipschitz
  continuous.}
\medskip

  To prove Theorem 2.1 we shall work in the path space
\begin{eqnarray*}
  \mathscr{X}_T &=&\{\mathbf{u}\in {L}^{\infty}_{\rm loc}((0,T),L^\infty(\mathbb{R}^N)):
  \nabla\cdot\mathbf{u}=0,\;\;
  \|\mathbf{u}\|_{\mathscr{X}_T}=\sup_{0<t<T}\sqrt{t}\Big|\ln\Big(\frac{t}{e^2T}\Big)
  \Big|\|\mathbf{u}(t)\|_{\infty}<\infty\}.
\end{eqnarray*}
  $(\mathscr{X}_T,\|\cdot\|_{\mathscr{X}_T})$ is a Banach space. The solution will lie
  in the space
\begin{eqnarray*}
  \mathscr{Y}_T=C_w([0,T],B^{-1(\ln)}_{\infty\infty}(\mathbb{R}^N))
  \cap\mathscr{X}_T.
\end{eqnarray*}
  To prove Theorem 2.2 we shall work in the path space
\begin{eqnarray*}
  \mathscr{X}_T^0 &=&\{\mathbf{u}\in\mathscr{X}_T:
  \lim_{t\to 0^+}\sqrt{t}\Big|\ln\Big(\frac{t}{e^2T}\Big)\Big|
  \|\mathbf{u}(t)\|_{\infty}=0\};
\end{eqnarray*}
  the solution will lie in the space
\begin{eqnarray*}
  \mathscr{Y}_T^0=C([0,T],\mathcal{C}^{-1(\ln)}_{bu}(\mathbb{R}^N))\cap\mathscr{X}_T^0.
\end{eqnarray*}
  We first prove some preliminary lemmas. In what follows, the expression
  $A\lesssim_T B$ means $A\leq C_TB$ for some constant $C_T$ depending on $T$; if
  $A\leq CB$ for some constant $C$ independent of $T$ then we write $A\lesssim B$.
\medskip

  {\bf Lemma 2.3}\ \ {\em If $\varphi\in B^{-1(\ln)}_{\infty\infty}(\mathbb{R}^N)$
  then $e^{t\Delta}\varphi\in\mathscr{Y}_T$ for any finite $T>0$, and
$$
  \|e^{t\Delta}\varphi\|_{\mathscr{X}_T}+\sup_{t\in (0,T)}
  \|e^{t\Delta}\varphi\|_{B^{-1(\ln)}_{\infty\infty}}\lesssim_T
  \|\varphi\|_{B^{-1(\ln)}_{\infty\infty}}.
$$
  If furthermore $\varphi\in\mathcal{C}^{-1(\ln)}_{bu}(\mathbb{R}^N)$
  then in addition to the above estimate we also have $e^{t\Delta}\varphi\in
  \mathscr{Y}_T^0$, and
$$
  \lim_{T\to 0^+}\|e^{t\Delta}\varphi\|_{\mathscr{X}_T}=0.
$$}

  {\em Proof}:\ \ It is easy to see that $B^{-1(\ln)}_{\infty\infty}(\mathbb{R}^N)$ is
  a shift-invariant Banach space of distributions. Hence by Propositions 4.1 and 4.4 of
  \cite{LEM02} we see that $\varphi\in B^{-1(\ln)}_{\infty\infty}(\mathbb{R}^N)$ implies
  that $e^{t\Delta}\varphi\in C_\ast([0,T],B^{-1(\ln)}_{\infty\infty}(\mathbb{R}^N))$, i.e.,
  the map $t\mapsto e^{t\Delta}\varphi$ from $[0,T]$ to $B^{-1(\ln)}_{\infty\infty}
  (\mathbb{R}^N)$ is continuous for $0<t\leq T$ with respect to the norm topology of
  $B^{-1(\ln)}_{\infty\infty}(\mathbb{R}^N)$ and continuous at $t=0$ with respect to the
  $\ast$-weak topology of $B^{-1(\ln)}_{\infty\infty}(\mathbb{R}^N)$, and
$$
  \sup_{t\in [0,T]}\|e^{t\Delta}\varphi\|_{B^{-1(\ln)}_{\infty\infty}}
  \leq\|\varphi\|_{B^{-1(\ln)}_{\infty\infty}}.
$$
  Moreover, from the definition of the space $B^{-1(\ln)}_{\infty\infty}(\mathbb{R}^N)$
  we see that $\varphi\in B^{-1(\ln)}_{\infty\infty}(\mathbb{R}^N)$ implies that
  $e^{t\Delta}\varphi\in\mathscr{X}_T$, and
$$
  \|e^{t\Delta}\varphi\|_{\mathscr{X}_T}\lesssim_T\|\varphi\|_{B^{-1(\ln)}_{\infty\infty}}.
$$
  Hence the first part of the lemma follows. Next we assume that
  $\varphi\in\mathcal{C}^{-1(\ln)}_{bu}(\mathbb{R}^N)$. We choose a sequence of functions
  $\psi_n\in C_{bu}(\mathbb{R}^N)$ ($n=1,2,\cdots$) such that $\lim_{n\to\infty}
  \|\varphi-\psi_n\|_{B^{-1(\ln)}_{\infty\infty}}=0$. We write
\begin{eqnarray*}
  \|e^{t\Delta}\varphi-\varphi\|_{B^{-1(\ln)}_{\infty\infty}}&\;\leq\;&
  \|e^{t\Delta}\varphi-e^{t\Delta}\psi_n\|_{B^{-1(\ln)}_{\infty\infty}}+
  \|e^{t\Delta}\psi_n-\psi_n\|_{B^{-1(\ln)}_{\infty\infty}}
  +\|\psi_n-\varphi\|_{B^{-1(\ln)}_{\infty\infty}}
\\
  &\leq& 2\|\varphi-\psi_n\|_{B^{-1(\ln)}_{\infty\infty}}
  +C\|e^{t\Delta}\psi_n-\psi_n\|_{\infty}.
\end{eqnarray*}
  Since $e^{t\Delta}$ is a $C_0$-semigroup in $C_{bu}(\mathbb{R}^N)$, we have $\lim_{t\to 0^+}
  \|e^{t\Delta}\psi_n-\psi_n\|_{\infty}=0$ (for each $n$). Hence, from the above estimate
  we immediately see that $\lim_{t\to 0^+}\|e^{t\Delta}\varphi-
  \varphi\|_{B^{-1(\ln)}_{\infty\infty}}=0$, so that $e^{t\Delta}\varphi$ is also
  continuous at $t=0$. Therefore, $e^{t\Delta}\varphi\in C([0,T],
  {\mathcal{C}}^{-1(\ln)}_{bu}(\mathbb{R}^N))$. Moreover, for any $\psi\in
  C_{bu}(\mathbb{R}^N)$ we can easily check that $e^{t\Delta}\psi\in\mathscr{X}_T^0$
  and $\lim_{T\to 0^+}\|e^{t\Delta}\psi\|_{\mathscr{X}_T}=0$. From this assertion and
  the fact that
$$
  \|e^{t\Delta}\varphi-e^{t\Delta}\psi_n\|_{\mathscr{X}_T}\lesssim
  \|\varphi-\psi_n\|_{B^{-1(\ln)}_{\infty\infty}}
$$
  (uniformly for $0<T\leq 1$) we see that also $e^{t\Delta}\varphi\in\mathscr{X}_T^0$
  and $\lim_{T\to 0^+}\|e^{t\Delta}\varphi\|_{\mathscr{X}_T}=0$. This completes the proof
  of the lemma. $\quad\Box$
\medskip

  In what follows we consider the bilinear map:
$$
  B(\mathbf{u},\mathbf{v})(t)=\int_0^te^{(t-\tau)\Delta}\mathbb{P}\nabla\cdot
  [\mathbf{u}(\tau)\otimes\mathbf{v}(\tau)]d\tau.
$$
  The following lemma will play a fundamental role:
\medskip

  {\bf Lemma 2.4} \ \ {\em Let $E$ be a shift-invariant Banach space of distributions.
  Let $\mathbf{u}$ be a $N\times N$ matrix function on $\mathbb{R}^N$ with all entries
  belonging to $E$. Then}
$$
  \|e^{t\Delta}\mathbb{P}\nabla\cdot\mathbf{u}\|_E\lesssim
  t^{-\frac{1}{2}}\|\mathbf{u}\|_E, \quad
  \mbox{for all}\;\; t>0.
$$

  {\em Proof}:\ \ In the case $E=L^p(\mathbb{R}^N)$ ($1\leq p\leq\infty$), this inequality
  was proved by Giga, Inui and Matsui in \cite{GIM99}. The general case follows from a
  similar argument. Indeed, by Proposition 11.1 of \cite{LEM02}, the kernel of the
  convolution operator $e^{t\Delta}\mathbb{P}\nabla$ belongs to $L^1(\mathbb{R}^N)$ and
  its $L^1$ norm is bounded by $C\sqrt{t}$. Hence the above inequality follows from the
  shift-invariance of $E$. $\quad\Box$
\medskip

  We now establish estimates on the bilinear map $B(\mathbf{u},\mathbf{v})$.
\medskip

  {\bf Lemma 2.5} \ \ {\em Let $T>0$ be given and assume that $\mathbf{u},\mathbf{v}
  \in{\mathscr{X}_T}$. Then $B(\mathbf{u},\mathbf{v})\in{\mathscr{X}_T}\cap
  C_w([0,T],B^{-1(\ln)}_{\infty\infty}(\mathbb{R}^N))$, and
\setcounter{equation}{0}
\begin{align}
  \|B(\mathbf{u},\mathbf{v})\|_{\mathscr{X}_T}
  &\lesssim\|\mathbf{u}\|_{\mathscr{X}_T}\|\mathbf{v}\|_{\mathscr{X}_T},
\\
  \sup_{t\in (0,T)}\|B(\mathbf{u},\mathbf{v})(t)&\|_{B^{-1(\ln)}_{\infty\infty}}
  \lesssim\|\mathbf{u}\|_{\mathscr{X}_T}\|\mathbf{v}\|_{\mathscr{X}_T}.
\\
  \lim_{t\to 0^+}B(\mathbf{u},\mathbf{v})(t)=0
  \;\; &(\mbox{$\ast$-weakly in}\;\;B^{-1(\ln)}_{\infty\infty}(\mathbb{R}^N)).
\end{align}
%--(2.1), (2.2),(2.3)--
  If furthermore either $\mathbf{u}\in\mathscr{Y}_T^0$ or $\mathbf{v}\in
  \mathscr{Y}_T^0$ then also $B(\mathbf{u},\mathbf{v})\in\mathscr{Y}_T^0$, and}
\begin{align}
  \lim_{t\to 0^+}\|B(\mathbf{u},\mathbf{v})(t)&\|_{B^{-1(\ln)}_{\infty\infty}}=0.
\end{align}
%--(2.4)--

  {\em Proof}:\ \ By Lemma 2.4 (choose $E=L^\infty(\mathbb{R}^N)$) we have
\begin{eqnarray}
  \|B(\mathbf{u},\mathbf{v})(t)\|_{\infty}
  &\lesssim&\int_0^t(t-\tau)^{-\frac{1}{2}}
  \|\mathbf{u}(\tau)\otimes\mathbf{v}(\tau)\|_{\infty}d\tau
\nonumber\\
  &\lesssim&\sup_{0<\tau<t}\sqrt{\tau}\Big|\ln\Big(\frac{\tau}{e^2T}\Big)\Big|
  \|\mathbf{u}(\tau)\|_{\infty}\cdot
  \sup_{0<\tau<t}\sqrt{\tau}\Big|\ln\Big(\frac{\tau}{e^2T}\Big)\Big|
  \|\mathbf{v}(\tau)\|_{\infty}
\nonumber\\
  &&\times\int_0^t(t-\tau)^{-\frac{1}{2}}\tau^{-1}
  \Big|\ln\Big(\frac{\tau}{e^2T}\Big)\Big|^{-2}d\tau.
%\nonumber\\
%  &\lesssim&\|\mathbf{u}\|_{\mathcal{L}^{\infty,\sigma}_T L_x^\infty}
%  \|\mathbf{v}\|_{\mathcal{L}^{\infty,\sigma}_T L_x^\infty}\cdot
%  \int_0^t(t-\tau)^{-\frac{1}{2}}\tau^{-1}
%  \Big|\ln\Big(\frac{\tau}{e^2T}\Big)\Big|^{-2}d\tau.
\end{eqnarray}
%--(2.5)--
  It is easy to show that
\begin{eqnarray}
  \int_0^t(t-\tau)^{-\frac{1}{2}}\tau^{-1}
  \Big|\ln\Big(\frac{\tau}{e^2T}\Big)\Big|^{-2}d\tau
  \lesssim  t^{-\frac{1}{2}}\Big|\ln\Big(\frac{t}{e^2T}\Big)\Big|^{-1}.
\end{eqnarray}
%--(2.6)--
  Indeed, we have
\begin{eqnarray*}
  &&\int_0^t(t-\tau)^{-\frac{1}{2}}\tau^{-1}
  \Big|\ln\Big(\frac{\tau}{e^2T}\Big)\Big|^{-2}d\tau
\\
  &=&\int_0^\frac{t}{2}(t-\tau)^{-\frac{1}{2}}\tau^{-1}
  \Big|\ln\Big(\frac{\tau}{e^2T}\Big)\Big|^{-2}d\tau
  +\int_\frac{t}{2}^t(t-\tau)^{-\frac{1}{2}}\tau^{-1}
  \Big|\ln\Big(\frac{\tau}{e^2T}\Big)\Big|^{-2}d\tau
\\
  &\lesssim&\sqrt{\frac{2}{t}}\int_0^{\frac{t}{2}}\tau^{-1}
  \Big|\ln\Big(\frac{\tau}{e^2T}\Big)\Big|^{-2}d\tau
  +\Big(\frac{t}{2}\Big)^{-1}\Big|\ln\Big(\frac{t}{2e^2T}\Big)\Big|^{-2}
  \int_{\frac{t}{2}}^t(t-\tau)^{-\frac{1}{2}}d\tau
\\
  &=&\Big(\frac{t}{2}\Big)^{-\frac{1}{2}}\Big|\ln\Big(\frac{t}{2e^2T}\Big)\Big|^{-1}
  +\sqrt{2t}\Big(\frac{t}{2}\Big)^{-1}\Big|\ln\Big(\frac{t}{2e^2T}\Big)\Big|^{-2}
\\
  &\lesssim&t^{-\frac{1}{2}}\Big|\ln\Big(\frac{t}{e^2T}\Big)\Big|^{-1}.
\end{eqnarray*}
%  Moreover, since $\lim_{t\to 0^+}\sqrt{t}|\ln(\frac{t}{e^2T})|
%  \|\mathbf{u}(t)\|_{\infty}=0$, from (2.4) and (2.5) we see that also\\
%  $\lim_{t\to 0^+}\sqrt{t}|\ln(\frac{t}{e^2T})|\|B(\mathbf{u},\mathbf{v})(t)\|_{\infty}=0$.
  Hence $B(\mathbf{u},\mathbf{v})(t)\in\mathscr{X}_T$ and (2.1) holds. Next,
  we note that also by Lemma 2.4, for any $s\in [0,T]$ and $t\in [0,T]$ we have
\begin{eqnarray*}
  \|e^{s\Delta}B(\mathbf{u},\mathbf{v})(t)\|_{\infty}
  &=&\Big\|\int_0^t\! e^{(t+s-\tau)\Delta}\mathbb{P}\nabla\cdot
  [\mathbf{u}(\tau)\otimes\mathbf{v}(\tau)]d\tau\Big\|_{\infty}
\\
  &\lesssim &\int_0^{t}\!(t+s-\tau)^{-\frac{1}{2}}
  \|\mathbf{u}(\tau)\otimes\mathbf{v}(\tau)\|_{\infty}d\tau
\\
  &\lesssim&\sup_{0<\tau<t}\sqrt{\tau}\Big|\ln\Big(\frac{\tau}{e^2T}\Big)\Big|
  \|\mathbf{u}(\tau)\|_{\infty}\cdot
  \sup_{0<\tau<t}\sqrt{\tau}\Big|\ln\Big(\frac{\tau}{e^2T}\Big)\Big|
  \|\mathbf{u}(\tau)\|_{\infty}
\nonumber\\
  &&\times\int_0^t(t+s-\tau)^{-\frac{1}{2}}\tau^{-1}
  \Big|\ln\Big(\frac{\tau}{e^2T}\Big)\Big|^{-2}d\tau.
\end{eqnarray*}
  By (2.6) we have
\begin{eqnarray*}
  \int_0^t(t+s-\tau)^{-\frac{1}{2}}\tau^{-1}
  \Big|\ln\Big(\frac{\tau}{e^2T}\Big)\Big|^{-2}d\tau
  &\lesssim &\int_0^{t+s}(t+s-\tau)^{-\frac{1}{2}}\tau^{-1}
  \Big|\ln\Big(\frac{\tau}{e^2T}\Big)\Big|^{-2}d\tau
\\
  &\lesssim & (t+s)^{-\frac{1}{2}}\Big|\ln\Big(\frac{t+s}{e^2T}\Big)\Big|^{-1}
\\
  &\lesssim & s^{-\frac{1}{2}}\Big|\ln\Big(\frac{s}{e^2T}\Big)\Big|^{-1}.
\end{eqnarray*}
  Hence for any $t\in [0,T]$ we have
\setcounter{equation}{6}
\begin{eqnarray}
  \|B(\mathbf{u},\mathbf{v})(t)\|_{B^{-1(\ln)}_{\infty\infty}}
  &=&\sup_{0<s<T}\sqrt{s}\Big|\ln\Big(\frac{s}{e^2T}\Big)\Big|
  \|e^{s\Delta}B(\mathbf{u},\mathbf{v})(t)\|_{\infty}
\nonumber\\
  &\lesssim &
  \sup_{0<\tau<t}\sqrt{\tau}\Big|\ln\Big(\frac{\tau}{e^2T}\Big)\Big|
  \|\mathbf{u}(\tau)\|_{\infty}\cdot
  \sup_{0<\tau<t}\sqrt{\tau}\Big|\ln\Big(\frac{\tau}{e^2T}\Big)\Big|
  \|\mathbf{v}(\tau)\|_{\infty}
\\
  &\lesssim &
  \|\mathbf{u}\|_{\mathscr{X}_T}\|\mathbf{v}\|_{\mathscr{X}_T}.
\nonumber
\end{eqnarray}
%---(2.7)---
  This proves (2.2).

  The proof of (2.3) is easy. Indeed, for any $\mathbf{w}\in\mathcal{S}(\mathbb{R}^N)$
  we have
\begin{eqnarray*}
  |\langle B(\mathbf{u},\mathbf{v})(t),\mathbf{w}\rangle|
  &=&|\int_0^t\langle e^{(t-\tau)\Delta}\mathbb{P}\nabla\cdot
  [\mathbf{u}(\tau)\otimes\mathbf{v}(\tau)],\mathbf{w}\rangle d\tau|
\\
  &\leq&\int_0^t\|\mathbf{u}(t)\|_{\infty}\|\mathbf{v}(t)\|_{\infty}
  \|\mathbb{P}\nabla\mathbf{w}\|_{1} d\tau
\\
  &\leq&\|\mathbf{u}\|_{\mathscr{X}_T}\|\mathbf{v}\|_{\mathscr{X}_T}
  \|\mathbb{P}\nabla\mathbf{w}\|_{1}\int_0^t\tau^{-1}
  \Big|\ln\Big(\frac{\tau}{e^2T}\Big)\Big|^{-2}d\tau
\\
  &=&\Big|\ln\Big(\frac{t}{e^2T}\Big)\Big|^{-1}
  \|\mathbf{u}\|_{\mathscr{X}_T}\|\mathbf{v}\|_{\mathscr{X}_T}
  \|\mathbb{P}\nabla\mathbf{w}\|_{1}\to 0 \;\;
  (\mbox{as}\;\; t\to 0^+).
\end{eqnarray*}
  Since $\mathcal{S}(\mathbb{R}^N)$ is dense in the predual of
  $B^{-1(\ln)}_{\infty\infty}(\mathbb{R}^N)$, using this result and (2.2) we
  obtain (2.3).

  We now proceed to prove that $B(\mathbf{u},\mathbf{v})\in C_w([0,T],
  B^{-1(\ln)}_{\infty\infty}(\mathbb{R}^N))$. By (2.3),
  $B(\mathbf{u},\mathbf{v})(t)$ is $\ast$-weakly continuous at $t=0$.
  We now consider an arbitrary point $0<t_0\leq T$. If $t_0<t<T$ then we write
\begin{eqnarray}
  &&B(\mathbf{u},\mathbf{v})(t)-B(\mathbf{u},\mathbf{v})(t_0)
\nonumber\\
  &=&\int_0^t e^{(t-\tau)\Delta}\mathbb{P}\nabla\cdot
  [\mathbf{u}(\tau)\otimes\mathbf{v}(\tau)]d\tau
  -\int_0^{t_0}e^{(t_0-\tau)\Delta}\mathbb{P}\nabla\cdot
  [\mathbf{u}(\tau)\otimes\mathbf{v}(\tau)]d\tau
\nonumber\\
  &=&\int_{t_0}^t e^{(t-\tau)\Delta}\mathbb{P}\nabla\cdot
  [\mathbf{u}(\tau)\otimes\mathbf{v}(\tau)]d\tau
  +[e^{(t-t_0)\Delta}-I]\int_0^{t_0}e^{(t_0-\tau)\Delta}\mathbb{P}\nabla\cdot
  [\mathbf{u}(\tau)\otimes\mathbf{v}(\tau)]d\tau
\nonumber\\
  &=:& A(t)+B(t),
\end{eqnarray}
%---(2.8)---
  and if $0<t_0-\delta<t<t_0$ then we write
\begin{eqnarray}
  &&B(\mathbf{u},\mathbf{v})(t_0)-B(\mathbf{u},\mathbf{v})(t)
\nonumber\\
  &=&\int_0^{t_0}e^{(t_0-\tau)\Delta}\mathbb{P}\nabla\cdot
  [\mathbf{u}(\tau)\otimes\mathbf{v}(\tau)]d\tau
  -\int_0^t e^{(t-\tau)\Delta}\mathbb{P}\nabla\cdot
  [\mathbf{u}(\tau)\otimes\mathbf{v}(\tau)]d\tau
\nonumber\\
  &=&\int_t^{t_0} e^{(t_0-\tau)\Delta}\mathbb{P}\nabla\cdot
  [\mathbf{u}(\tau)\otimes\mathbf{v}(\tau)]d\tau
  +[e^{(t_0-t)\Delta}-I]\int_0^t e^{(t-\tau)\Delta}\mathbb{P}\nabla\cdot
  [\mathbf{u}(\tau)\otimes\mathbf{v}(\tau)]d\tau
\nonumber\\
  &=&\int_t^{t_0} e^{(t_0-\tau)\Delta}\mathbb{P}\nabla\cdot
  [\mathbf{u}(\tau)\otimes\mathbf{v}(\tau)]d\tau
  +[e^{(t_0-t)\Delta}-I]\int_{t_0-\delta}^t e^{(t-\tau)\Delta}\mathbb{P}\nabla\cdot
  [\mathbf{u}(\tau)\otimes\mathbf{v}(\tau)]d\tau
\nonumber\\
  &&+e^{(t-t_0+\delta)\Delta}[e^{(t_0-t)\Delta}-I]\int_0^{t_0-\delta}
  e^{(t_0-\delta-\tau)\Delta}\mathbb{P}\nabla\cdot
  [\mathbf{u}(\tau)\otimes\mathbf{v}(\tau)]d\tau
\nonumber\\
  &=:& A_1(t)+B_1(t)+B_2(t).
\end{eqnarray}
%---(2.9)---
  For $A(t)$ we have (see the proof of (2.7))
\begin{eqnarray*}
  \|A(t)\|_{B^{-1(\ln)}_{\infty\infty}}&=&\sup_{0<s<T}\sqrt{s}
  \Big|\ln\Big(\frac{s}{e^2T}\Big)\Big|
  \|\int_{t_0}^t e^{(t+s-\tau)\Delta}\mathbb{P}\nabla\cdot
  [\mathbf{u}(\tau)\otimes\mathbf{v}(\tau)]d\tau\|_{\infty}
\\
  &\lesssim &
  \|\mathbf{u}\|_{\mathscr{X}_T}\|\mathbf{v}\|_{\mathscr{X}_T}
  \sup_{0<s<T}\sqrt{s}\Big|\ln\Big(\frac{s}{e^2T}\Big)\Big|
  \int_{t_0}^t(t+s-\tau)^{-\frac{1}{2}}\tau^{-1}
  \Big|\ln\Big(\frac{\tau}{e^2T}\Big)\Big|^{-2}d\tau
\\
  &\lesssim_T&
  \|\mathbf{u}\|_{\mathscr{X}_T}\|\mathbf{v}\|_{\mathscr{X}_T}
  \int_{t_0}^t(t-\tau)^{-\frac{1}{2}}\tau^{-1}
  \Big|\ln\Big(\frac{\tau}{e^2T}\Big)\Big|^{-2}d\tau,
\end{eqnarray*}
  so that $\displaystyle\lim_{t\to t_0^+}\|A(t)\|_{B^{-1(\ln)}_{\infty\infty}}=0$.
  Moreover, since $B(t)=[e^{(t-t_0)\Delta}-I]B(\mathbf{u},\mathbf{v})(t_0)$ and
  $B(\mathbf{u},\mathbf{v})(t_0)\in B^{-1(\ln)}_{\infty\infty}(\mathbb{R}^N)$, by Lemma
  2.3 we have $\displaystyle\lim_{t\to t_0^+}B(t)=0$ ($\ast$-weakly in
  $B^{-1(\ln)}_{\infty\infty}(\mathbb{R}^N)$). Hence
$$
  \lim_{t\to t_0^+}B(\mathbf{u},\mathbf{v})(t)=B(\mathbf{u},\mathbf{v})(t_0)
  \;\; (\mbox{$\ast$-weakly in}\;\;B^{-1(\ln)}_{\infty\infty}(\mathbb{R}^N)).
$$
  Next, similarly as for $A(t)$ we have $\displaystyle\lim_{t\to t_0^-}
  \|A_1(t)\|_{B^{-1(\ln)}_{\infty\infty}}=0$. Moreover, similarly as for the treatment
  of $A(t)$ we have that by choosing $\delta$ sufficiently small,
  $\|B_1(t)\|_{B^{-1(\ln)}_{\infty\infty}}$ can be as small as we expect, and when
  $\delta$ is chosen and fixed, $B_2(t)$ can be treated similarly as for $B(t)$ to
  get that for any $\mathbf{w}$ in the predual of $B^{-1(\ln)}_{\infty\infty}
  (\mathbb{R}^N)$, $\displaystyle\lim_{t\to t_0^-}\langle B_2(t),\mathbf{w}\rangle=0$.
  Hence
$$
  \lim_{t\to t_0^-}B(\mathbf{u},\mathbf{v})(t)=B(\mathbf{u},\mathbf{v})(t_0)
  \;\; (\mbox{$\ast$-weakly in}\;\;B^{-1(\ln)}_{\infty\infty}(\mathbb{R}^N)).
$$
  This proves $B(\mathbf{u},\mathbf{v})\in C_w([0,T],B^{-1(\ln)}_{\infty\infty}
  (\mathbb{R}^N))$.

  We now furthermore assume that either $\mathbf{u}\in\mathscr{Y}_T^0$ or $\mathbf{v}
  \in\mathscr{Y}_T^0$. It is not difficult to see that in this case we have
  $B(\mathbf{u},\mathbf{v})(t)\in\mathcal{C}^{-1(\ln)}_{bu}(\mathbb{R}^N)$ for each
  fixed $t\in [0,T]$, and from (2.7) we see that (2.4) holds. Moreover, since
  $e^{t\Delta}$ is a $C_0$-semigroup in $C_{bu}(\mathbb{R}^N)$ and $B(\mathbf{u},
  \mathbf{v})(t_0)\in\mathcal{C}^{-1(\ln)}_{bu}(\mathbb{R}^N)$ for any $0<t_0\leq T$,
  from the definition of $B(t)$ (see (2.8)) we see that $\lim_{t\to t_0^+}
  \|B(t)\|_{B^{-1(\ln)}_{\infty\infty}}=0$, so that
$$
  \lim_{t\to t_0^+}\|B(\mathbf{u},\mathbf{v})(t)-B(\mathbf{u},
  \mathbf{v})(t_0)\|_{B^{-1(\ln)}_{\infty\infty}}=0
$$
  for any $0<t_0<T$. Similarly we have
$$
  \lim_{t\to t_0^-}\|B(\mathbf{u},\mathbf{v})(t)-
  B(\mathbf{u},\mathbf{v})(t_0)\|_{B^{-1(\ln)}_{\infty\infty}}=0.
$$
  for any $0<t_0\leq T$. Hence $B(\mathbf{u},\mathbf{v})\in C([0,T],
  \mathcal{C}^{-1(\ln)}_{bu}(\mathbb{R}^N))$. Finally, from (2.5) and (2.6) we
  see that $B(\mathbf{u},\mathbf{v})\in\mathscr{X}_T^0$. It follows that $B(\mathbf{u},
  \mathbf{v})\in\mathscr{Y}_T^0$. This completes the proof of Lemma 2.6. $\quad\Box$
\medskip

  Having proved Lemmas 2.3 and 2.5, Theorems 2.1 and 2.2 follow from a standard fixed
  point argument. Indeed, we rewrite the problem (1.3) into the following equivalent
  integral equation:
$$
  \mathbf{u}(t)=e^{t\Delta}\mathbf{u}_0-B(\mathbf{u},\mathbf{u})(t).
$$
  We define a map $\mathscr{F}$ such that for given $\mathbf{u}$, $\mathscr{F}(\mathbf{u})$
  equals to the right-hand side of the above equation. By Lemmas 2.3 and 2.5, if
  $\mathbf{u}_0\in B^{-1(\ln)}_{\infty\infty}(\mathbb{R}^N)$ then $\mathscr{F}$ maps
  $\mathscr{X}_T$ into $\mathscr{Y}_T\subseteq\mathscr{X}_T$. Furthermore, we have
  the following estimates:
\begin{eqnarray*}
  \|\mathscr{F}(\mathbf{u})\|_{\mathscr{X}_T}&\,\leq\,&
  \|e^{t\Delta}\mathbf{u}_0\|_{\mathscr{X}_T}+C\|\mathbf{u}\|_{\mathscr{X}_T}^2,
\\
  \|\mathscr{F}(\mathbf{u})-\mathscr{F}(\mathbf{v})\|_{\mathscr{X}_T}&\,\leq\,&
  C(\|\mathbf{u}\|_{\mathscr{X}_T}+\|\mathbf{v}\|_{\mathscr{X}_T})
  \|\mathbf{u}-\mathbf{v}\|_{\mathscr{X}_T}.
\end{eqnarray*}
  To prove the assertion $(i)$ of Theorem 2.1, for given $T>0$ we let $\mathbf{u}_0\in
  B^{-1}_{\infty 2}(\mathbb{R}^N)$ (with $\diva\mathbf{u}_0=0$) be so small that
  $\|e^{t\Delta}\mathbf{u}_0\|_{\mathscr{X}_T}\leq\epsln$,
  where $\epsln$ is an arbitrarily chosen positive number such that $4C\epsln<1$. Then
  from the first inequality we easily see that $\mathscr{F}$ maps the closed ball
  $\overline{B}(0,2\epsln)$ in $\mathscr{X}_T$ into itself, and the second inequality
  ensures that $\mathscr{F}$ is a contraction mapping when restricted to this ball. Hence,
  by the fixed point theorem of Banach, $\mathscr{F}$ has a unique fixed point in this
  ball. Since $e^{t\Delta}\mathbf{u}_0\in\mathscr{Y}_T$ and $\mathscr{F}(\mathscr{X}_T)
  \subseteq\mathscr{Y}_T$, we get a mild solution of the problem (1.1) which lies in the
  space $\mathscr{Y}_T$. This proves Theorem 2.1. The proof of Theorem 2.2 is similar,
  and we omit it here.

\section{Linearly perturbed Navier-Stokes equations}

\hskip 2em
  We now consider the following initial value problem:
\setcounter{equation}{0}
\begin{eqnarray}
\left\{
\begin{array}{l}
  \partial_t\mathbf{u}-\Delta\mathbf{u}+\mathbb{P}\nabla\cdot(2\,\mathbf{a}\otimes_s
  \mathbf{u}+\mathbf{u}\otimes\mathbf{u})=0,\quad x\in\mathbb{R}^N,\;\; t>0,\\
  \mbox{div}\,\mathbf{u}=0,\qquad x\in\mathbb{R}^N,\;\; t>0,\\
\mathbf{u}(x,0)=\mathbf{u}_{0}(x),\quad x\in\mathbb{R}^N,
\end{array}
\right.
\end{eqnarray}
%--(3.1)--
  where $\mathbf{a}=\mathbf{a}(x,t)$ is a given vector function in $\mathbb{R}^N\times
  \mathbb{R}_+$, and $\mathbf{a}\otimes_s\mathbf{u}=\displaystyle\frac{1}{2}
  (\mathbf{a}\otimes\mathbf{u}+\mathbf{u}\otimes\mathbf{a})$.

  In this section we study well-posedness of the above problem. We shall not consider
  the most general situation, but for the purpose of our later application we only
  consider the case where the initial data belongs to $B_{\widetilde{X}_r}^{s_r,q_r}
  (\mathbb{R}^N)$ for some $0<r<1$,  where $s_r=-1+r$ and $q_r=\frac{2}{1-r}$.
\medskip

  {\bf Theorem 3.1}\ \ {\em Let $0<T<\infty$, $\mathbf{a}$ and $\mathbf{u}_0$ be given such
  that $\sqrt{t}\mathbf{a}(t)\in L^\infty(\mathbb{R}^N\times[0,T])$ and $\mathbf{u}_0\in
  B_{{X}_r}^{s_r,q_r}(\mathbb{R}^N)$, where $0<r<1$, $s_r=-1+r$ and $q_r=\frac{2}{1-r}$.
  Assume that $\displaystyle\lim_{t\to 0^+}\sqrt{t}\|\mathbf{a}(t)\|_{\infty}=0$.
  Then we have the following assertions:

  $(i)$ Given $\mathbf{a}$ and $\mathbf{u}_0$ as above, there exists corresponding $T_1
  =T_1(\mathbf{a},\mathbf{u}_0)\in (0,T]$ such that the problem $(3.1)$ has a unique
  mild solution in the class
\begin{eqnarray*}
  &\mathbf{u}\in C([0,T_1],B_{{X}_r}^{s_r,q_r}(\mathbb{R}^N))
  \cap L^{q_r}([0,T_1],{X}_r(\mathbb{R}^N)), \quad
  \nabla\cdot\mathbf{u}=0, &
\\
  &t^{\frac{1-r}{2}}\mathbf{u}(t)\in L^\infty([0,T_1],{X}_r(\mathbb{R}^N)), \quad
  \sqrt{t}\mathbf{u}(t)\in L^\infty([0,T_1],L^\infty(\mathbb{R}^N)),&
\\
  &\displaystyle\lim_{t\to 0^+}t^{\frac{1-r}{2}}\|\mathbf{u}(t)\|_{{X}_r}=0, \quad
  \lim_{t\to 0^+}\sqrt{t}\|\mathbf{u}(t)\|_{\infty}=0.&
\end{eqnarray*}
  Moreover, the solution map $\mathbf{u}_0\mapsto\mathbf{u}$ from a small neighborhood
  of any given point $\mathbf{v}_0$ of $B_{{X}_r}^{s_r,q_r}(\mathbb{R}^N)$ to the Banach
  space of functions of the above class $($with $T_1=T_1(\mathbf{a},\mathbf{v}_0))$ is
  Lipschitz continuous. Besides, if $\mathbf{u}_0\in B_{\widetilde{X}_r}^{s_r,q_r}
  (\mathbb{R}^N)$ then $\mathbf{u}\in C([0,T_1],B_{\widetilde{X}_r}^{s_r,q_r}
  (\mathbb{R}^N))$.

  $(ii)$ Given $T>0$, there exist corresponding constants $\epsln_1>0$ and $\epsln_2>0$
  depending only on $T$ and $N$, such that if $\sup_{t\in [0,T]}\!\sqrt{t}
  \|\mathbf{a}(t)\|_{\infty}<\epsln_1$ and $\|\mathbf{u}_0\|_{B_{{X}_r}^{s_r,q_r}}<
  \epsln_2$, the problem $(3.1)$ has a unique mild solution in the class
\begin{eqnarray*}
  &\mathbf{u}\in C([0,T],B_{{X}_r}^{s_r,q_r}(\mathbb{R}^N))
  \cap L^{q_r}([0,T],{X}_r(\mathbb{R}^N)), \quad
  \nabla\cdot\mathbf{u}=0, &
\\
  &t^{\frac{1-r}{2}}\mathbf{u}(t)\in L^\infty([0,T],{X}_r(\mathbb{R}^N)), \quad
  \sqrt{t}\mathbf{u}(t)\in L^\infty([0,T],L^\infty(\mathbb{R}^N)),&
\\
  &\displaystyle\lim_{t\to 0^+}t^{\frac{1-r}{2}}\|\mathbf{u}(t)\|_{{X}_r}=0, \quad
  \lim_{t\to 0^+}\sqrt{t}\|\mathbf{u}(t)\|_{\infty}=0.&
\end{eqnarray*}
  Moreover, the solution map $\mathbf{u}_0\mapsto\mathbf{u}$ from the ball $\overline{B}
  (0,\epsln)$ of $B_{{X}_r}^{s_r,q_r}(\mathbb{R}^N)$ to the Banach space of functions of
  the above class is Lipschitz continuous. Besides, if $\mathbf{u}_0\in
  B_{\widetilde{X}_r}^{s_r,q_r}(\mathbb{R}^N)$ then $\mathbf{u}\in
  C([0,T],B_{\widetilde{X}_r}^{s_r,q_r}(\mathbb{R}^N))$.}
\medskip

  To prove this result, we rewrite the problem (3.1) into the following equivalent
  integral equation:
$$
  \mathbf{u}(t)=e^{t\Delta}\mathbf{u}_0-B(\mathbf{a},\mathbf{u})
  -B(\mathbf{u},\mathbf{a})-B(\mathbf{u},\mathbf{u}).
\eqno{(3.2)}
$$
  To prove existence of a solution to this equation, we shall work in the path space
\begin{eqnarray*}
  \mathscr{V}_T=\Big\{\mathbf{u}\in L^{q_r}([0,T],&&{X}_r(\mathbb{R}^N)):
   \sup_{t\in [0,T]}t^{\frac{1-r}{2}}\|u(t)\|_{{X}_r}<\infty,
  \sup_{t\in [0,T]}\sqrt{t}\|u(t)\|_{\infty}<\infty,
\\
  &&\;\;\lim_{t\to 0^+}t^{\frac{1-r}{2}}\|\mathbf{u}(t)\|_{{X}_r}=0, \;
  \lim_{t\to 0^+}\sqrt{t}\|\mathbf{u}(t)\|_{\infty}=0\Big\},
\end{eqnarray*}
  which is a Banach space with respect to the norm
$$
  \|u\|_{\mathscr{V}_T}=\|u\|_{L_T^{q_r}{X}_r}+
  \sup_{t\in [0,T]}t^{\frac{1-r}{2}}\|u(t)\|_{{X}_r}+
  \sup_{t\in [0,T]}\sqrt{t}\|u(t)\|_{\infty},
$$
  where $\|u\|_{L_T^{q_r}{X}_r}=\displaystyle\Big(\int_0^T\!\|u(t)\|_{{X}_r}^{q_r}dt
  \Big)^{\frac{1}{q_r}}$. We shall prove that the solution lies in the space
%  We shall prove that $\mathscr{G}$ maps $\mathscr{V}_T$ into itself, and is
%  a contraction mapping when restricted to some closed ball in $\mathscr{V}_T$, so that
%  it has a fixed point in $\mathscr{V}_T$. In fact, we shall prove that the image of
%  $\mathscr{G}$ is contained in the space
\begin{eqnarray*}
  \mathscr{W}_T=C([0,T],B_{{X}_r}^{s_r,q_r}(\mathbb{R}^N))
  \cap\mathscr{V}_T.
\end{eqnarray*}

  As in Section 2, the proof will be based on several preliminary Lemmas.
\medskip

  {\bf Lemma 3.2}\ \ {\em $\{e^{t\Delta}\}_{t\geq 0}$ is a $C_0$-semigroup in
  $B_{{X}_r}^{s_r,q_r}(\mathbb{R}^N)$, and for any $\varphi\in
  B_{{X}_r}^{s_r,q_r}(\mathbb{R}^N)$ we have $e^{t\Delta}\varphi\in\mathscr{W}_T$, and
$$
  \|e^{t\Delta}\varphi\|_{\mathscr{W}_T}\lesssim_T\|\varphi\|_{B_{{X}_r}^{s_r,q_r}}.
$$
  Moreover, if $\varphi\in B_{\widetilde{X}_r}^{s_r,q_r}(\mathbb{R}^N)$ then $e^{t\Delta}
  \varphi\in C([0,T],B_{\widetilde{X}_r}^{s_r,q_r}(\mathbb{R}^N))$.}
\medskip

  {\em Proof}:\ \ Since $B_{{X}_r}^{s_r,q_r}(\mathbb{R}^N)$ is a shift-invariant Banach
  space of distributions, for any $\varphi\in B_{{X}_r}^{s_r,q_r}(\mathbb{R}^N)$ we have
$$
  \|e^{t\Delta}\varphi\|_{B_{{X}_r}^{s_r,q_r}}\lesssim_T
  \|\varphi\|_{B_{{X}_r}^{s_r,q_r}},
$$
  and $e^{t\Delta}\varphi\in C_\ast([0,T],B_{{X}_r}^{s_r,q_r}(\mathbb{R}^N))$. Moreover,
  the heat kernel characterization of $B_{{X}_r}^{s_r,q_r}(\mathbb{R}^N)$ implies
  that
$$
  \|e^{t\Delta}\varphi\|_{L^{q_r}_T{X}_r}\lesssim_T\|\varphi\|_{B_{{X}_r}^{s_r,q_r}},
$$
  the embedding $B_{{X}_r}^{s_r,q_r}(\mathbb{R}^N)\subseteq B_{{X}_r}^{s_r,\infty}
  (\mathbb{R}^N)$ and the heat kernel characterization of $B_{{X}_r}^{s_r,\infty}
  (\mathbb{R}^N)$ implies that
$$
  \sup_{t\in [0,T]}t^{\frac{1-r}{2}}\|e^{t\Delta}\varphi\|_{{X}_r}
  \lesssim_T\|\varphi\|_{B_{{X}_r}^{s_r,q_r}},
$$
  and the embedding $B_{{X}_r}^{s_r,q_r}(\mathbb{R}^N)\hookrightarrow B^{-1}_{\infty\infty}
  (\mathbb{R}^N)$ and the heat kernel characterization of $B^{-1}_{\infty\infty}
  (\mathbb{R}^N)$ implies that
$$
  \sup_{t\in [0,T]}\sqrt{t}\|e^{t\Delta}\varphi\|_{\infty}\lesssim
  \|\varphi\|_{B_{{X}_r}^{s_r,q_r}}.
$$
  The embedding $B_{{X}_r}^{s_r,q_r}(\mathbb{R}^N)\hookrightarrow B^{-1}_{\infty\infty}
  (\mathbb{R}^N)$ follows from the embedding ${X}_r(\mathbb{R}^N)\hookrightarrow
  B^{-r}_{\infty\infty}(\mathbb{R}^N)$, whose proof can be found in \cite{Ger06}.

  Next we prove that for any $\varphi\in B_{{X}_r}^{s_r,q_r}(\mathbb{R}^N)$,
  $\displaystyle\lim_{t\to 0^+}\|e^{t\Delta}\varphi-\varphi\|_{B_{{X}_r}^{s_r,q_r}}=0$.
  Indeed, since $q_r<\infty$, the Littlewood-Paley decomposition of $\varphi$ implies
  that $\cap_{m=0}^\infty H^m_{{X}_r}(\mathbb{R}^N)$ is dense in $B_{{X}_r}^{s_r,q_r}
  (\mathbb{R}^N)$, where $H^m_{{X}_r}(\mathbb{R}^N)=(I-\Delta)^{-\frac{m}{2}}{X}_r
  (\mathbb{R}^N)$. Since ${X}_r(\mathbb{R}^N)\hookrightarrow B^{-r}_{\infty\infty}
  (\mathbb{R}^N)$, we have $H^m_{{X}_r}(\mathbb{R}^N)\hookrightarrow B^{m,\infty}_{{X}_r}
  (\mathbb{R}^N)\hookrightarrow B^{m-r}_{\infty\infty}(\mathbb{R}^N)\hookrightarrow
  L^{\infty}(\mathbb{R}^N)$ for $m\geq 1$, so that $H^{m+1}_{{X}_r}(\mathbb{R}^N)
  \hookrightarrow W^{m,\infty}(\mathbb{R}^N)$ for $m\geq 1$. Thus $C_{bu}(\mathbb{R}^N)$ is
  dense in $B_{{X}_r}^{s_r,q_r}(\mathbb{R}^N)$. Now for any $\varphi\in
  B_{{X}_r}^{s_r,q_r}(\mathbb{R}^N)$ and $\psi\in C_{bu}(\mathbb{R}^N)\subseteq
  L^{\infty}(\mathbb{R}^N)\subseteq {X}_r(\mathbb{R}^N)$, we write
\begin{eqnarray*}
  \|e^{t\Delta}\varphi-\varphi\|_{B_{{X}_r}^{s_r,q_r}}&\;\leq\;&
  \|e^{t\Delta}\varphi-e^{t\Delta}\psi\|_{B_{{X}_r}^{s_r,q_r}}+
  \|e^{t\Delta}\psi-\psi\|_{B_{{X}_r}^{s_r,q_r}}+\|\psi-\varphi\|_{B_{{X}_r}^{s_r,q_r}}
\\
  &\leq& 2\|\varphi-\psi\|_{B_{{X}_r}^{s_r,q_r}}+\|e^{t\Delta}\psi-\psi\|_{{X}_r}
\\
  &\leq& 2\|\varphi-\psi\|_{B_{{X}_r}^{s_r,q_r}}+\|e^{t\Delta}\psi-\psi\|_{\infty}.
\end{eqnarray*}
  From this inequality and a similar argument as in the proof of Lemma 2.3, we see that
  the desired assertion follows. Finally, by the density of $L^{\infty}(\mathbb{R}^N)$
  and ${X}_r(\mathbb{R}^N)$ in $B_{{X}_r}^{s_r,q_r}(\mathbb{R}^N)$ and a similar argument
  as in the proof of Lemma 2.3, we easily deduce that for any $\varphi\in
  B_{{X}_r}^{s_r,q_r}(\mathbb{R}^N)$,
$$
  \lim_{t\to 0^+}t^{\frac{1-r}{2}}\|e^{t\Delta}\varphi\|_{{X}_r}=0
  \quad \mbox{and} \quad
  \lim_{t\to 0^+}\sqrt{t}\|e^{t\Delta}\varphi\|_{\infty}=0.
$$
  This proves the lemma. $\quad\Box$
\medskip

  The last three terms on the right-hand side of (3.2) will be treated with the same type
  of estimates. To establish such estimates, we need the following delicate inequality
  (see Lemma 20.1 of \cite{LEM02}):
\medskip

  {\bf Lemma 3.3} \ \ {\em Let $\theta\in (0,\frac{1}{2})$. Then for any $1\leq q\leq\infty$,
  the operator $f\mapsto g$ defined by
$$
  g(t)=\int_0^t\frac{1}{\sqrt{(t-\tau)\tau}}\Big(\frac{t}{\tau}\Big)^{\theta}f(\tau)d\tau
$$
  is bounded on $L^q((0,T),\frac{dt}{t})$. More precisely, there exists a constant $C>0$
  depending only on $\theta$ and $q$ $($independent of $f$ and $T)$ such that
$$
  \Big(\int_0^T|g(t)|^q\frac{dt}{t}\Big)^{\frac{1}{q}}\leq
  C\Big(\int_0^T|f(t)|^q\frac{dt}{t}\Big)^{\frac{1}{q}}\quad
  \mbox{if}\;\; 1\leq q<\infty,
$$
  and $\|g\|_{\infty}\leq C\|f\|_{\infty}$ if $q=\infty$. $\quad\Box$}
\medskip

  {\bf Lemma 3.4} \ \ {\em Assume that $\sqrt{t}\mathbf{v}(t)\in L^\infty([0,T],
  L^\infty(\mathbb{R}^N))$. Then we have the following assertions:

  $(i)$ If $\mathbf{u}\in L^{q_r}([0,T],{X}_r(\mathbb{R}^N))$ then $B(\mathbf{u},
  \mathbf{v})$, $B(\mathbf{v},\mathbf{u})\in L^{q_r}([0,T],{X}_r(\mathbb{R}^N))\cap
  C([0,T],B_{{X}_r}^{s_r,q_r}(\mathbb{R}^N))$, and
\setcounter{equation}{2}
\begin{eqnarray}
  \|B(\mathbf{u},\mathbf{v})\|_{L_T^{q_r}{X}_r}+
  \|B(\mathbf{v},\mathbf{u})\|_{L_T^{q_r}{X}_r}
  &\;\lesssim\;&
  \sup_{t\in [0,T]}\!\sqrt{t}\|\mathbf{v}(t)\|_{\infty}
  \|\mathbf{u}\|_{L_T^{q_r}{X}_r},
\\
  \sup_{t\in [0,T]}\|B(\mathbf{u},\mathbf{v})(t)\|_{B_{{X}_r}^{s_r,q_r}}
  +\sup_{t\in [0,T]}\|B(\mathbf{v},\mathbf{u})\|_{B_{{X}_r}^{s_r,q_r}}
  &\;\lesssim\;&
  \sup_{t\in [0,T]}\!\sqrt{t}\|\mathbf{v}(t)\|_{\infty}
  \|\mathbf{u}\|_{L_T^{q_r} {X}_r}.
\end{eqnarray}
%--(3.4)--

  $(ii)$ If $t^{\frac{1-r}{2}}\mathbf{u}(t)\in L^\infty([0,T],
  {X}_r(\mathbb{R}^N))$ then $t^{\frac{1-r}{2}}B(\mathbf{u},\mathbf{v})$,
  $t^{\frac{1-r}{2}}B(\mathbf{v},\mathbf{u})\in L^\infty([0,T],{X}_r(\mathbb{R}^N))$,
  $\sqrt{t}B(\mathbf{u},\mathbf{v})$, $\sqrt{t}B(\mathbf{v},\mathbf{u})\in
  L^\infty([0,T],L^\infty(\mathbb{R}^N))$, and}
\begin{eqnarray}
  \sup_{t\in [0,T]}\!t^{\frac{1-r}{2}}\|B(\mathbf{u},\mathbf{v})\|_{{X}_r}\!+\!
  \sup_{t\in [0,T]}\!t^{\frac{1-r}{2}}\|B(\mathbf{v},\mathbf{u})\|_{{X}_r}
  &\;\lesssim\;&
  \sup_{t\in [0,T]}\!\sqrt{t}\|\mathbf{v}(t)\|_{\infty}
  \sup_{t\in [0,T]}\!t^{\frac{1-r}{2}}\|\mathbf{u}(t)\|_{{X}_r},
\\
  \sup_{t\in [0,T]}\!\sqrt{t}\|B(\mathbf{u},\mathbf{v})(t)\|_{\infty}\!+\!
  \sup_{t\in [0,T]}\!\sqrt{t}\|B(\mathbf{v},\mathbf{u})(t)\|_{\infty}
  &\;\lesssim\;&
  \sup_{t\in [0,T]}\!\sqrt{t}\|\mathbf{v}(t)\|_{\infty}
  \sup_{t\in [0,T]}\!t^{\frac{1-r}{2}}\|\mathbf{u}(t)\|_{{X}_r}.\;\;\;
\end{eqnarray}
%--(3.5), (3.6)--

  {\em Proof}:\ \ The proofs of (3.3), (3.5) and (3.6) are contained in the proof of
  Theorem 20.2 of \cite{LEM02}. Here we only give the proof of (3.4) and the assertion
  that $B(\mathbf{u},\mathbf{v})$, $B(\mathbf{v},\mathbf{u})\in C([0,T],
  B_{{X}_r}^{s_r,q_r}(\mathbb{R}^N))$. Moreover, we only give the
  proof for $B(\mathbf{v},\mathbf{u})$, because that for $B(\mathbf{u},\mathbf{v})$
  is similar. We first prove (3.4). We use the heat kernel characterization of
  $B_{{X}_r}^{s_r,q_r}(\mathbb{R}^N)$. First we have
\begin{eqnarray*}
  \|e^{\Delta}B(\mathbf{v},\mathbf{u})(t)\|_{{X}_r}
  &\lesssim&\int_0^t(1+t-\tau)^{-\frac{1}{2}}\|\mathbf{v}(\tau)\|_{\infty}
  \|\mathbf{u}(\tau)\|_{{X}_r}d\tau
\\
  &\lesssim&\sup_{t\in [0,T]}\!\sqrt{t}\|\mathbf{v}(t)\|_{\infty}
  \int_0^t(1+t-\tau)^{-\frac{1}{2}}\tau^{-\frac{1}{2}}
  \|\mathbf{u}(\tau)\|_{{X}_r}d\tau
\\
  &\lesssim&\sup_{t\in [0,T]}\!\sqrt{t}\|\mathbf{v}(t)\|_{\infty}
  \cdot\Big(\int_0^\infty|1-\tau|^{-\frac{q_r'}{2}}\tau^{-\frac{q_r'}{2}}
  d\tau\Big)^{\frac{1}{q_r'}}\|\mathbf{u}\|_{L_T^{q_r} {X}_r}
\\
  &\lesssim&\sup_{t\in [0,T]}\!\sqrt{t}\|\mathbf{v}(t)\|_{\infty}
  \|\mathbf{u}\|_{L_T^{q_r} {X}_r}, \quad \forall t\in [0,T].
\end{eqnarray*}
  Next, let $\mathbf{u}^\ast$ be as before. Then for any $0<s<1$ and
  $0\leq t\leq T$ we have
\begin{eqnarray*}
  \|e^{s\Delta}B(\mathbf{v},\mathbf{u})(t)\|_{{X}_r}
  &\lesssim&\int_0^t(t+s-\tau)^{-\frac{1}{2}}\|\mathbf{v}(\tau)\|_{\infty}
  \|\mathbf{u}(\tau)\|_{{X}_r}d\tau
\\
  &\lesssim&\sup_{t\in [0,T]}\!\sqrt{t}\|\mathbf{v}(t)\|_{\infty}
  \int_0^{t+s}(t+s-\tau)^{-\frac{1}{2}}\tau^{-\frac{1}{2}}
  \|\mathbf{u}^\ast(\tau)\|_{{X}_r}d\tau,
\end{eqnarray*}
  so that
\begin{eqnarray*}
  s^{\frac{1-r}{2}}\|e^{s\Delta}B(\mathbf{v},\mathbf{u})(t)\|_{{X}_r}
  &\lesssim&\sup_{t\in [0,T]}\!\sqrt{t}\|\mathbf{v}(t)\|_{\infty}
  \int_0^{t+s}(t+s-\tau)^{-\frac{1}{2}}\tau^{-\frac{1}{2}}
  \Big(\frac{t+s}{\tau}\Big)^{\frac{1-r}{2}}\cdot
  \tau^{\frac{1-r}{2}}\|\mathbf{u}^\ast(\tau)\|_{{X}_r}d\tau.
\end{eqnarray*}
  Using Lemma 3.3, we get
\begin{eqnarray*}
  \Big[\int_0^1\|e^{s\Delta}B(\mathbf{v},\mathbf{u})(t)\|_{{X}_r}^{q_r}ds
  \Big]^{\frac{1}{q_r}}
  &=&\Big[\int_0^1\Big(s^{\frac{1-r}{2}}\|e^{s\Delta}
  B(\mathbf{v},\mathbf{u})(t)\|_{{X}_r}\Big)^{q_r}\frac{ds}{s}
  \Big]^{\frac{1}{q_r}}
\\
  &\lesssim&\sup_{t\in [0,T]}\!\sqrt{t}\|\mathbf{v}(t)\|_{\infty}
  \Big[\int_0^\infty\Big(\tau^{\frac{1-r}{2}}
  \|\mathbf{u}^\ast(\tau)\|_{{X}_r}\Big)^{q_r}\frac{d\tau}{\tau}
  \Big]^{\frac{1}{q_r}}
\\
  &\lesssim&\sup_{t\in [0,T]}\!\sqrt{t}\|\mathbf{v}(t)\|_{\infty}
  \|\mathbf{u}\|_{L_T^{q_r} {X}_r}.
\end{eqnarray*}
  This proves (3.4).

  Having proved (3.4), the assertion that $B(\mathbf{v},\mathbf{u})\in
  C([0,T],B_{{X}_r}^{s_r,q_r}(\mathbb{R}^N))$ follows from a similar argument
  as in the proof of the corresponding assertion in Lemma 2.6. $\quad\Box$
\medskip

  We are now ready to give the proof of Theorem 3.1.

  First, by Lemma 3.4 we see that for a given function $\mathbf{a}$ on $\mathbb{R}^N\times
  [0,T]$ such that $\sqrt{t}\mathbf{a}(t)\in L^\infty([0,T],L^\infty(\mathbb{R}^N))$, the
  map
$$
  L(\mathbf{a}):\mathbf{u}\mapsto-B(\mathbf{a},\mathbf{u})-B(\mathbf{u},\mathbf{a})
$$
  is a bounded linear operator in $\mathscr{V}_T$, with norm bounded by $C\sup_{t\in
  [0,T]}\!\sqrt{t}\|\mathbf{a}(t)\|_{\infty}$, where $C$ is a positive constant
  independent of $\mathbf{a}$ and $T$. It follows that if
$$
  C\sup_{t\in [0,T]}\!\sqrt{t}\|\mathbf{a}(t)\|_{\infty}\leq\frac{1}{2},
\eqno{(3.7)}
$$
  then $I-L(\mathbf{a})$ is invertible, and $\|[I-L(\mathbf{a})]^{-1}\|_{\mathcal{L}
  (\mathscr{V}_T)}\leq 2$. Clearly, the condition (3.7) is satisfied in each of the
  following two situations:
\begin{itemize}
\item Given $\mathbf{a}$, restrict our discussion in a possibly smaller interval
  $[0,T_1]$ ($0<T_1\leq T$) so that\\ $C\sup_{t\in [0,T_1]}\!\sqrt{t}\|\mathbf{a}(t)
  \|_{\infty}\leq\frac{1}{2}$. Existence of a such $T_1$ is guaranteed by the condition
  $\displaystyle\lim_{t\to 0^+}\sqrt{t}\|\mathbf{a}(t)\|_{\infty}=0$. In this case
  we need to replace all $T$ with $T_1$ in the succeeding discussion.
\item Given $T>0$, restrict our discussion to those $\mathbf{a}$ such that
  $C\sup_{t\in [0,T_1]}\!\sqrt{t}\|\mathbf{a}(t)\|_{\infty}\leq\frac{1}{2}$.
  This is a smallness assumption on $\mathbf{a}$.
\end{itemize}
  In both situations, it can be easily seen that when considering solutions in
  $\mathscr{V}_T$, the equation (3.2) is equivalent to the following equation:
$$
  \mathbf{u}(t)=[I-L(\mathbf{a})]^{-1}[e^{t\Delta}\mathbf{u}_0-B(\mathbf{u},\mathbf{u})].
\eqno{(3.8)}
$$
  Letting $\mathscr{G}(\mathbf{u})$ be the right-hand side of (3.8), by using Lemmas 3.2
  and 3.4 we see that $\mathscr{G}$ maps $\mathscr{V}_T$ to itself, and the following
  estimates hold:
\begin{eqnarray*}
  \|\mathscr{G}(\mathbf{u})\|_{\mathscr{V}_T}&\,\leq\,&
  2\|e^{t\Delta}\mathbf{u}_0\|_{\mathscr{V}_T}+C\|\mathbf{u}\|_{\mathscr{V}_T}^2,
\\
  \|\mathscr{G}(\mathbf{u})-\mathscr{G}(\mathbf{v})\|_{\mathscr{V}_T}&\,\leq\,&
  C[\|\mathbf{u}\|_{\mathscr{V}_T}+\|\mathbf{v}\|_{\mathscr{V}_T}]
  \|\mathbf{u}-\mathbf{v}\|_{\mathscr{V}_T}.
\end{eqnarray*}
  Thus, a similar argument as in the proof of Theorem 2.1 shows that $\mathscr{G}$ has a
  fixed point in a small closed neighborhood of the origin of $\mathscr{V}_T$ in each of
  the following two situations:
\begin{itemize}
\item Given $\mathbf{u}_0\in B_{{X}_r}^{s_r,q_r}(\mathbb{R}^N)$, restrict
  our discussion in a possibly smaller interval $[0,T_2]$ ($0<T_2\leq T_1$) so that
  $2\|e^{t\Delta}\mathbf{u}_0\|_{\mathscr{V}_{T_2}}\leq\epsln$, where $\epsln$ is an
  arbitrarily chosen positive number such that $4C\epsln<1$. Existence of a such number
  $T_2$ is ensured by Lemma 3.2. In this case we need to replace all $T$ with $T_2$ in
  the succeeding discussion.
\item Given $T>0$, restrict our discussion to those $\mathbf{u}_0\in
  B_{{X}_r}^{s_r,q_r}(\mathbb{R}^N)$ such that $2C_0\|\mathbf{u}_0\|_{B_{{X}_r}^{s_r,q_r}}
  \leq\epsln$, where $C_0$ is the constant that appears in the inequality $\|e^{t\Delta}
  \mathbf{u}_0\|_{\mathscr{V}_T}\leq C_0\|\mathbf{u}_0\|_{B_{{X}_r}^{s_r,q_r}}$, and
  $\epsln$ is as above. This is a smallness assumption on $\mathbf{u}_0$.
\end{itemize}
  Consequently, we obtain a solution of (3.2) in $\mathscr{V}_T$. Since by
  Lemmas 3.2 and 3.4, for any $\mathbf{u}_0\in B_{{X}_r}^{s_r,q_r}(\mathbb{R}^N)$ and
  $\mathbf{u}\in\mathscr{V}_T$, all terms on the right-hand side of (3.2) belong to
  $\mathscr{W}_T$, we actually have obtained a solution of (3.2) in $\mathscr{W}_T$.
  Moreover, since $e^{t\Delta}$ is also a $C_0$-semigroup in $B_{\widetilde{X}_r}^{s_r,q_r}
  (\mathbb{R}^N)$, a similar discussion shows that if $\mathbf{u}_0\in
  B_{\widetilde{X}_r}^{s_r,q_r}(\mathbb{R}^N)$ then the solution $\mathbf{u}$ belongs to
  $C([0,T],B_{\widetilde{X}_r}^{s_r,q_r}(\mathbb{R}^N))\cap\mathscr{V}_T$. This completes
  the proof of Theorem 3.1.  $\quad\Box$

\section{The proof of Theorem 1.4}

\hskip 2em
  In this section we give the proof of Theorem 1.4. We first prove the following result:
\medskip

  {\bf Theorem 4.1}\ \ {\em Let $0<T<\infty$ and $\mathbf{a}=\mathbf{a}_1+\mathbf{a}_2$
  be given, where $\mathbf{a}_1\in L^2([0,T],L^\infty(\mathbb{R}^N))$ and $\mathbf{a}_2\in
  L^{q_r}([0,T],{X}_r(\mathbb{R}^N))$ $(0<r<1$, $q_r=\frac{2}{1-r})$. Then for any
  $\mathbf{u}_0\in L^2(\mathbb{R}^N)$ satisfying the condition $\nabla\cdot\mathbf{u}_0=0$,
  the problem $(3.1)$ has a weak solution $\mathbf{u}\in C_w([0,T],L^2(\mathbb{R}^N))\cap
  L^2([0,T],{H}^1(\mathbb{R}^N))$ satisfying the following generalized energy inequality:
\setcounter{equation}{0}
\begin{eqnarray}
  \sup_{0\leqslant t\leqslant T}\|\mathbf{u}(t)\|_2^2+\int_0^T\!\!
  \|\nabla\mathbf{u}(t)\|_2^2dt\leq\|\mathbf{u}_0\|_2^2
  \exp\Big[C\!\int_0^t\!\!\Big(\|\mathbf{a}_1(\tau)\|_{\infty}^2
  +\|\mathbf{a}_2(\tau)\|_{{X}_r}^{\frac{2}{1-r}}\Big)d\tau\Big].
\end{eqnarray}
%---(4.1)---
}

  {\em Proof}:\ \ Choose a sequence of divergence-free vector functions $\{\mathbf{w}_k
  \}_{k=1}^\infty\subseteq\cap_{m=0}^\infty H^m(\mathbb{R}^N)$ such that they form a
  normalized orthogonal basis of $\mathbb{P}L^2(\mathbb{R}^N)$. Let $\mathbf{u}_{0}=
  \sum_{k=0}^\infty c_k\mathbf{w}_k$ be the Fourier expansion of $\mathbf{u}_{0}$ with respect
  to this basis of $\mathbb{P}L^2(\mathbb{R}^N)$. For each $n$ we define an approximate
  solution $\mathbf{u}_n$ of the problem (3.1) on the time interval $[0,T]$ as follows:
$$
  \mathbf{u}_n(x,t)=\sum_{k=1}^n g_{kn}(t)\mathbf{w}_k(x), \quad x\in\mathbb{R}^N, \quad t\in [0,T],
$$
  where $\{g_{kn}\}_{k=1}^n$ is the solution of the initial value problem
\begin{eqnarray}
\left\{
\begin{array}{l}
  \displaystyle g_{kn}'(t)=\sum_{j=1}^n b_{jk}(t)g_{jn}(t)+\sum_{i,j=1}^n c_{ijk}g_{in}(t)g_{jn}(t),
  \quad 0<t<T, \quad k=1,2,\cdots,n,
\\
  \displaystyle g_{kn}(0)=c_k, \quad k=1,2,\cdots,n,
\end{array}
\right.
\end{eqnarray}
%--(4.2)--
  where
\begin{eqnarray*}
\left\{
\begin{array}{l}
  \displaystyle b_{jk}(t)=-\int_{\mathbb{R}^N}\nabla\mathbf{w}_j(x)
  \cdot\nabla\mathbf{w}_k(x)dx+
  2\int_{\mathbb{R}^N}[\mathbf{a}(x,t)\otimes_s\mathbf{w}_j(x)]
  \cdot[\nabla\otimes\mathbf{w}_k(x)]dx,
\\
  \displaystyle c_{ijk}=\int_{\mathbb{R}^N}[\mathbf{w}_i(x)\otimes\mathbf{w}_j(x)]
  \cdot[\nabla\otimes\mathbf{w}_k(x)]dx
\end{array}
\right.
\end{eqnarray*}
  ($i,j,k=1,2,\cdots,n$). Note that by the assumption on $\mathbf{a}$ we have $b_{jk}\in
  L^2(0,T)$ ($i,j,k=1,2,\cdots,n$), so that the problem (4.2) has a unique solution
  $\{g_{kn}\}_{k=1}^n$ defined in a maximal interval $I$ which is either $[0,T]$ or
  $[0,T^\ast)$ with $0<T^\ast\leq T$, such that $g_{kn}\in C(I)$ and $g_{kn}'\in
  L^2_{{\footnotesize loc}}(I)$, $k=1,2,\cdots,n$, and in the second case $g_{kn}(t)$
  blows up as $t\to T^{\ast-}$. By making a priori estimates we can ensure that the second
  case cannot occur. Indeed, (4.2) can be rewritten as follows (note that $\mathbb{P}
  \Delta\mathbf{u}_n=\Delta\mathbb{P}\mathbf{u}_n=\Delta\mathbf{u}_n$):
\begin{eqnarray}
\left\{
\begin{array}{l}
  \partial_t\mathbf{u}_n-\Delta\mathbf{u}_n+\mathbb{P}_n[2\nabla\cdot(\mathbf{a}
  \otimes_s\mathbf{u}_n)+\nabla\cdot(\mathbf{u}_n\otimes\mathbf{u}_n)]=0,
  \quad x\in\mathbb{R}^N,\;\; t>0,\\
  \mathbf{u}_n(x,0)=\mathbf{u}_{0n}(x),\quad x\in\mathbb{R}^N,
\end{array}
\right.
\end{eqnarray}
%--(4.3)--
  where $\mathbb{P}_n$ denotes the projection from $L^2(\mathbb{R}^N)$ to its $n$-dimensional
  subspace spanned by $\{\mathbf{w}_k\}_{k=1}^n$, and $\mathbf{u}_{0n}=\sum_{k=1}^n c_k\mathbf{w}_k$.
  It is easy to see that
\begin{eqnarray*}
  \int_{\mathbb{R}^N}[\mathbf{a}_1(x,t)\otimes_s\mathbf{u}_n(x,t)]\cdot
  [\nabla\otimes\mathbf{u}_n(x,t)]dx
  &\;\leq\; &\|\mathbf{a}_1(t)\|_{\infty}\|\mathbf{u}_n(t)\|_2\|\nabla\mathbf{u}_n(t)\|_2,
\\
  \int_{\mathbb{R}^N}[\mathbf{a}_2(x,t)\otimes_s\mathbf{u}_n(x,t)]\cdot
  [\nabla\otimes\mathbf{u}_n(x,t)]dx
%  &\leq &\|\mathbf{a}_2(t)\|_{X_r}\|\mathbf{u}_n(t)\|_{H^r}\|\nabla\mathbf{u}_n(t)\|_2
%\\
  &\;\leq\; &\|\mathbf{a}_2(t)\|_{\dot{X}_r}\|\mathbf{u}_n(t)\|_2^{1-r}
  \|\nabla\mathbf{u}_n(t)\|_2^{1+r}.
\end{eqnarray*}
  Using these inequalities and a standard argument we can easily prove that
\begin{eqnarray*}
  \frac{1}{2}\frac{d}{dt}\|\mathbf{u}_n(t)\|_2^2
  &\;\leq\; &-\frac{1}{2}\|\nabla\mathbf{u}_n(t)\|_2^2+C\Big(\|\mathbf{a}_1(t)\|_{\infty}^2
  +\|\mathbf{a}_2(t)\|_{\dot{X}_r}^{\frac{2}{1-r}}\Big)\|\mathbf{u}_n(t)\|_2^2.
\end{eqnarray*}
  Hence
\begin{eqnarray}
  \|\mathbf{u}_n(t)\|_2^2+\int_0^t\!\!
  \|\nabla\mathbf{u}_n(\tau)\|_2^2d\tau\leq\|\mathbf{u}_{0n}\|_2^2
  \exp\Big[C\!\int_0^t\!\!\Big(\|\mathbf{a}_1(\tau)\|_{\infty}^2
  +\|\mathbf{a}_2(\tau)\|_{\dot{X}_r}^{\frac{2}{1-r}}\Big)d\tau\Big].
\end{eqnarray}
%--(4.4)--
  From this estimate, the desired assertion follows immediately.

  By (4.4), we easily deduce that for any $2\leq p\leq\frac{2N}{N-2}$ and $2\leq q\leq
  \infty$ satisfying the relation $\frac{N}{p}+\frac{2}{q}=\frac{N}{2}$,
  $\{\mathbf{u}_n\}_{n=1}^\infty$ is bounded in $L^q([0,T],L^p(\mathbb{R}^N))\cap
  L^2([0,T],H^1(\mathbb{R}^N))$. This implies that $\{\nabla\cdot(\mathbf{a}\otimes_s
  \mathbf{u}_n)\}_{n=1}^\infty$, $\{\nabla\cdot(\mathbf{u}_n\otimes
  \mathbf{u}_n)\}_{n=1}^\infty$ and $\{\Delta\mathbf{u}_n\}_{n=1}^\infty$ are bounded in
  $L^2([0,T],H^{-2}(\mathbb{R}^N))$. Using the equation (4.3), we further deduce that
  $\{\partial_t\mathbf{u}_n\}_{n=1}^\infty$ is bounded in $L^2([0,T],H^{-2}(\mathbb{R}^N))$.
  From these assertions and Theorem 2.1 in Chapter III of \cite{TEM}, we conclude that there
  is a subsequence $\{\mathbf{u}_{n_k}\}_{k=1}^\infty$ of $\{\mathbf{u}_n\}_{n=1}^\infty$,
  and a divergence-free vector function $\mathbf{u}$ such that
\begin{itemize}
\item $\mathbf{u}_{n_k}\to\mathbf{u}$ $*$-weakly in $L^\infty([0,T],L^2(\mathbb{R}^N))$ and
  weakly in $L^2([0,T],H^1(\mathbb{R}^N))$;
\item $\partial_t\mathbf{u}_{n_k}\to\partial_t\mathbf{u}$ weakly in $L^2([0,T],H^{-2}(\mathbb{R}^N))$;
\item For any bounded measurable set $Q\subseteq\mathbb{R}^N$ and any $2\leq p<\frac{2N}{N-2}$,
  $\mathbf{u}_{n_k}\to\mathbf{u}$ strongly in $L^2([0,T],L^p(Q))$.
\end{itemize}
  These assertions imply that $\mathbf{u}\in L^\infty([0,T],L^2(\mathbb{R}^N))\cap
  L^2([0,T],H^1(\mathbb{R}^N))$ and $\partial_t\mathbf{u}\in L^2([0,T],H^{-2}(\mathbb{R}^N))$,
  which further implies that $\mathbf{u}\in C([0,T],H^{-1/2}(\mathbb{R}^N))$, so that
  $\mathbf{u}\in C_w([0,T],L^2(\mathbb{R}^N))$. Moreover, by using the last assertion listed
  above, the boundedness of $\{\mathbf{u}_n\}_{n=1}^\infty$ in $L^{\frac{2(N+2)}{N}}
  (\mathbb{R}^N\times[0,T])$, and the Vitali convergence theorem (cf. Corollary A.2 of
  \cite{CK}), we conclude that for any $1\leq p<2(N\!+\!2)/N$ and any bounded measurable
  set $Q\subseteq\mathbb{R}^N$, $\mathbf{u}_{n_k}\to\mathbf{u}$ strongly in $L^p(Q\times
  [0,T])$. This immediately implies that for any $1\leqslant p<(N\!+\!2)/N$ and any bounded
  measurable set $Q\subseteq\mathbb{R}^N$, $\mathbf{u}_{n_k}\otimes\mathbf{u}_{n_k}\to
  \mathbf{u}\otimes\mathbf{u}$ strongly in $L^p(Q\times [0,T])$. It follows that
\begin{eqnarray*}
  \mathbf{u}_{n_k}\otimes\mathbf{u}_{n_k}\to\mathbf{u}\otimes\mathbf{u} \;\;
  \mbox{weakly in distribution sense},
\end{eqnarray*}
  so that also $\mathbf{u}_{n_k}\otimes\mathbf{u}_{n_k}\to\mathbf{u}\otimes\mathbf{u}$
  weakly in $L^{1+\frac{2}{N}}(\mathbb{R}^N\times[0,T])$,
%\begin{eqnarray*}
%  \mathbf{u}_{n_k}\otimes\mathbf{u}_{n_k}\to\mathbf{u}\otimes\mathbf{u} \;\;
%  \mbox{weakly in} \;\; L^{\frac{N+2}{N}}(\mathbb{R}^N\times[0,T])
%\end{eqnarray*}
  because $\mathbf{u}_{n_k}\otimes\mathbf{u}_{n_k}$ is bounded in $L^{1+\frac{2}{N}}
  (\mathbb{R}^N\times[0,T])$.
  Having obtained this assertion, we can now let $n=n_k\to\infty$ in (4.3) to conclude that
  $\mathbf{u}$ is a weak solution of the problem (3.1) in $\mathbb{R}^N\times [0,T]$. The
  inequality (4.1) is an immediate consequence of (4.4) and the Fatou lemma. The proof of
  Theorem 3.1 is complete. $\quad\Box$
\medskip

  {\em Proof of Theorem 1.4}:\ \ Let $0<T<\infty$ be given. We first assume that $\mathbf{u}_0
  =\mathbf{v}_0+\mathbf{w}_0+\mathbf{z}_0$, $\mathbf{v}_0\in{B}^{-1(\ln)}_{\infty\infty}
  (\mathbb{R}^N)$, $\mathbf{w}_0\in{B}_{\dot{X}_r}^{s_r,q_r}(\mathbb{R}^N)$, $\mathbf{z}_0
  \in L^2(\mathbb{R}^N)$, $\nabla\cdot\mathbf{v}_0=\nabla\cdot\mathbf{w}_0=\nabla\cdot
  \mathbf{z}_0=0$, $\|\mathbf{v}_0\|_{{B}^{-1(\ln)}_{\infty\infty}}<\epsln_1$ and
  $\|\mathbf{w}_0\|_{{B}_{{X}_r}^{s_r,q_r}}<\epsln_2$, where $\epsln_1$ and $\epsln_2$
  are positive numbers to be specified later. Consider the problem (1.1) with $\mathbf{u}_0$
  replaced by $\mathbf{v}_0$. By Theorem 2.1, we see that if $\epsln_1$ is sufficiently small
  then that problem has a unique solution which we denote as $\mathbf{v}$, such that $\mathbf{v}
  \in C([0,T], {B}^{-1(\ln)}_{\infty\infty}(\mathbb{R}^N))\cap\mathscr{X}_T$. Note that
  $\mathbf{v}\in\mathscr{X}_T$ implies that $\mathbf{v}\in L^2([0,T],L^\infty(\mathbb{R}^N))$,
  $\sqrt{t}\mathbf{v}\in L^\infty(\mathbb{R}^N\times[0,T])$ and $\displaystyle
  \lim_{t\to 0^+}\sqrt{t}\|\mathbf{v}(t)\|_{\infty}=0$. Moreover, from the proof of
  Theorem 2.1 we see that $\sup_{0\leq t\leq T}\|\mathbf{v}(t)\|_{{B}^{-1(\ln)}_{\infty\infty}}+
  \sup_{0\leq t\leq T}\sqrt{t}\|\mathbf{v}(t)\|_{\infty}+\|\mathbf{v}\|_{L_T^2 L_x^\infty}\leq
  C\epsln_1$, where $C$ is a constant depending only on the dimension $N$. Next consider
  the problem (3.1) with $\mathbf{a}$ and $\mathbf{u}_0$ replaced by $\mathbf{v}$
  and $\mathbf{w}_0$, respectively. By Theorem 3.1 $(ii)$, we see that if $\epsln_1$ and
  $\epsln_2$ are sufficiently small then that problem has a unique solution which
  we denote as $\mathbf{w}$, such that $\mathbf{w}\in C([0,T],{B}_{{X}_r}^{s_r,q_r}
  (\mathbb{R}^N))\cap L^{q_r}([0,T],{X}_r(\mathbb{R}^N))$ (we neglect the other properties
  of $\mathbf{w}$ because we do not use them below). We now consider the problem (3.1) again,
  but replace $\mathbf{a}$ and $\mathbf{u}_0$ with $\mathbf{v}+\mathbf{w}$ and $\mathbf{z}_0$,
  respectively. By using Theorem 4.1 with $\mathbf{a}_1=\mathbf{v}$ and $\mathbf{a}_2=
  \mathbf{w}$, we see that that problem has a weak solution, which we denote as
  $\mathbf{z}$, such that $\mathbf{z}\in C_w([0,T],L^2(\mathbb{R}^N))\cap
  L^2([0,T],H^1(\mathbb{R}^N))$. Now let $\mathbf{u}=\mathbf{v}+\mathbf{w}+
  \mathbf{z}$. Then as one can easily verify, $\mathbf{u}$ is the desired solution of the
  original problem (1.1). This proves the first assertion of Theorem 1.4.

  The second assertion is an easy consequence of the first one. Indeed, since
  $C_0^\infty(\mathbb{R}^N)$ is dense in $\mathcal{C}^{-1(\ln)}_{0}(\mathbb{R}^N)$, by
  using Proposition 12.1 of \cite{LEM02} we can
  split $\mathbf{v}_0$ into a sum $\mathbf{v}_0=\mathbf{v}_0'+\mathbf{v}_0''$, such that
  $\mathbf{v}_0'\in\mathcal{C}^{-1(\ln)}_{0}(\mathbb{R}^N)$, $\mathbf{v}_0''\in
  L^2(\mathbb{R}^N)$, $\nabla\cdot\mathbf{v}_0'=\nabla\cdot\mathbf{v}_0''=0$, and
  $\|\mathbf{v}_0'\|_{{B}^{-1(\ln)}_{\infty\infty}}$ is as small as we want. Similarly,
  since $C_0^\infty(\mathbb{R}^N)$ is dense in ${B}_{\widetilde{{X}}_r}^{s_r,q_r}
  (\mathbb{R}^N)$, we can split $\mathbf{w}_0$ into a sum $\mathbf{w}_0=\mathbf{w}_0'+
  \mathbf{w}_0''$, such that $\mathbf{w}_0'\in{B}_{\widetilde{{X}}_r}^{s_r,q_r}
  (\mathbb{R}^N)$, $\mathbf{w}_0''\in L^2(\mathbb{R}^N)$, $\nabla\cdot\mathbf{w}_0'=
  \nabla\cdot\mathbf{w}_0''=0$, and $\|\mathbf{w}_0'\|_{{B}_{\widetilde{{X}}_r}^{s_r,q_r}}$
  is as small as we want. Thus, by using the first assertion we immediately
  obtain the second assertion. This completes the proof of Theorem 1.4. $\quad\Box$

\end{document}